\documentclass[12pt]{amsart}
\usepackage{amssymb}
\usepackage{amsmath}
\usepackage{mathabx}
\usepackage{amsthm}
\usepackage{amstext}
\usepackage{amsopn}
\usepackage{mathrsfs}
\usepackage{latexsym}
\usepackage{enumerate}
\usepackage{bm}
\DeclareMathAlphabet{\mathpzc}{OT1}{pzc}{m}{it}
\allowdisplaybreaks
\newtheorem{Definition}{Definition}[section]
\newtheorem{Proposition}{Proposition}[section]
\newtheorem{Lemma}{Lemma}[section]
\newtheorem{Theorem}{Theorem}[section]
\newtheorem{Corollary}{Corollary}[section]
\newtheorem{Remark}{Remark}[section]
\newtheorem{Example}{Example}[section]

\begin{document}
\bibliographystyle{plain}
\footnotetext{
This work is supported by Narodowe Centrum Nauki grant No. 2014/15/B/ST1/00166}
\title[Motzkin path decompositions of functionals]
{Motzkin path decompositions of functionals\\in noncommutative probability}
\author[R. Lenczewski]{Romuald Lenczewski}
\address{Romuald Lenczewski \newline
Katedra Matematyki, Politechnika Wroc\l{}awska, \newline
Wybrze\.{z}e Wyspia\'{n}skiego 27, 50-370 Wroc{\l}aw, Poland}
\email{Romuald.Lenczewski@pwr.edu.pl}
\begin{abstract}
We study the decomposition of free random variables in terms of their {\it ortho\-go\-nal replicas} from a new perspective.
First, we show that the mixed moments of orthogonal replicas with respect to the normalized linear functional $\Phi$ 
are naturally described in terms of Motzkin paths identified with reduced Motzkin words $\mathpzc{M}$. Using this fact, we demonstrate that the mixed moments of order $n$ of free random variables with respect to the free product of normalized linear functionals 
are sums of the mixed moments of the orthogonal replicas of these variables with respect to $\Phi$ with 
summation extending over $\mathpzc{M}_{n}$, the set of reduced Motzkin paths of lenght $n$.
One of the applications of this formula is a decomposition formula for mixed moments of 
free random variables in terms of their boolean cumulants which corresponds to the decomposition of
the lattice ${\rm NC}(n)$ into sublattices $\mathcal{M}(w)$ of partitions which are monotonically adapted to 
colors in $w\in \mathpzc{M}_{n}$.
The linear functionals defined by the mixed moments of orthogonal replicas and indexed by reduced Motzkin words
play the role of a generating set of the space of product functionals in which the boolean product corresponds 
to constant Motzkin paths and the free product corresponds to all
Motzkin paths.
\end{abstract}

\maketitle

\section{Introduction}
Different notions of independence in noncommutative probability lead to different theories which 
are usually studied separately. In our unified approach to independence, in which
{\it free probability} of Voiculescu \cite{[25],[26],[27]} plays a special role, we are 
looking for a nice algebraic framework that would allow us to look at different notions of independence and the associated objects simultaneously. In this paper, we initiate a study based on the lattices of Motzkin paths \cite{[M]}.
In particular, we show that it is natural to consider families of linear as well as multilinear functionals indexed by Motzkin paths.
Our previous results on the {\it unification} of independence \cite{[14],[15],[16]}
and on the {\it decomposition} of free random variables \cite{[17],[18]} are closely related to this work, although 
the present approach seems to put them in a rather new perspective. 

The main objects of interest in noncommutative probability theories 
are random variables and their distributions, 
studied in terms of moments and cumulants related to different notions 
of independence. 
We constructed in \cite{[14],[16]} a unified tensor 
product noncommutative probability space in which different random variables were included. 
In particular, {\it boolean independent} random variables were the first order approximations of 
free random variables expressed in terms of series of {\it orthogonal replicas}. 
Moreover, random variables associated with other interesting 
notions of noncommmutative independence were also included in this scheme.
It is remarkable that all of them are in one way or another related to freeness.
Thus, two types of `independent' random variables which appear in the 
context of the decomposition of the free additive convolution and 
are important for the subordination property in free probability, 
namely {\it orthogonally independent} and {\it free with subordination} (or, {\it s-free}) random variables,
were expressed in \cite{[17]} in terms of orthogonal replicas.
These two notions are intrinsic to the decompositions of the free product of graphs, as observed in \cite{[1]} and \cite{[17]}, 
as well as to the recent operad approach to independence of Jekel and Liu \cite{[10]}.
Finally, in the case of two algebras, {\it monotone independent} random variables
can also be expressed in terms of orthogonal replicas. 
This follows from our recent paper \cite {[L2019]} on the more general 
independence of Hasebe called conditionally monotone \cite{[H]}. In the general case, 
a modification of the present approach is required and it will be treated elsewhere.

Therefore, the scheme of orthogonal replicas allows us to study moments of various types of variables simultaneously. 
In this paper, we implement it to introduce a common algebraic framework for a family of 
multilinear (and associated linear) moment functionals. In the simplest case of two algebras, 
${\mathcal A}_{1}$ and ${\mathcal A}_{2}$, the main idea relies on the decomposition of free random variables of Voiculescu of the form
$$
A=\sum_{j=1}^{\infty}a(j)\;\;\;{\rm and}\;\;\;B=\sum_{j=1}^{\infty}b(j)
$$
and the sequences 
$$
\{a(j):j\in {\mathbb N}\}\;\;\;{\rm and}\;\;\; \{b(j):j\in {\mathbb N}\}
$$ 
consist of suitably constructed replicas of variables 
$a\in {\mathcal A}_{1}$, $b\in {\mathcal A}_{2}$ called {\it orthogonal replicas} (or, simply {\it replicas}), 
respectively. 
We can treat these series as sequences of partial sums that approximate free random variables, 
which is sufficient for computations of all moments \cite{[14]}.
In the case of *-nonco\-m\-mu\-ta\-tive probability spaces one can also attribute a meaning to the series
on the algebraic level, using the concept of operators with monotone closure, similar to
Berberian's operators with closure (OWC) for Baer *-rings \cite{[5]}, as we showed in \cite{[16]}. 
In the $C^*$-algebra context, we get strongly convergent series of bounded operators on a Hilbert space.
Here, we detour this issue by treating functionals rather than variables.

We show that this decomposition leads to Motzkin paths in a quite natural way. 
In the case of two algebras, the mixed moments of random variables 
which are free with respect to some normalized linear functional are completely described in terms
of mixed moments of variables from $\mathcal{A}_{1}$ and $\mathcal{A}_{2}$ with respect to a family  
$$
\{\psi(w):w\in \mathpzc{M}^{*}\}
$$ 
of path-dependent linear moment functionals, called {\it Motzkin linear functionals}, defined on the 
free product ${\mathcal A}_{1}*{\mathcal A}_{2}$ without identification of units, where
$\mathpzc{M}^{*}:=\mathpzc{M}\cup \{\varnothing\}$ and $\mathpzc{M}$ is the family
of Motzkin paths, $\varnothing$ denoting the empty path.
These functionals are defined by a certain {\it duality relation} in terms of mixed 
moments of orthogonal replicas:
$$
\psi(w)\left(a_1\cdots a_n\right):=\Phi\left(a_1(j_1)\cdots a_n(j_n)\right)
$$
with respect to some tensor product state $\Phi$, where $a_k\in \mathcal{A}_{i_k}$,
for any $k\in [n]:=\{1, \ldots, n\}$, where $i_1\neq \cdots \neq i_n$ and
$$
w:=j_1\cdots j_n
$$ 
is a {\it reduced Motzkin word}, by which we understand a word in letters from the alphabet $\mathbb{N}$
such that $j_1=j_n=1$ and $|j_k-j_{k-1}|\in\{0,1\}$ for $1<k \leq n$. Reduced Motzkin words are
in bijection with Motzkin paths and for that reason we will denote both sets by $\mathpzc{M}$.

Moreover, if $K_1\subset {\mathcal A}_{1}$ and $K_2\subset {\mathcal A}_{2}$, then 
sets $\{a(1):a\in K_1\}\;\;\;{\rm and}\;\;\;\{b(1):y\in K_2\}$
are boolean independent with respect to $\Phi$. This leads us to the conclusion that the family of 
functionals corresponding to constant Motzkin paths
$$
\{\psi(1^{n}):n\in \mathbb{N}^{*}\}
$$
describes the boolean product of functionals, where $\mathbb{N}^{*}=\mathbb{N}\cup\{0\}$.
More importantly, the free product of normalized functionals 
is reproduced only if we take functionals corresponding to all Motzkin paths.
In fact, this is one of the main results of this paper, stated in the 
functional language in Theorem 5.1.

Using the language of mixed moments, it says that the mixed moments of 
free random variables under the free product of normalized linear functionals 
$\star_{i\in I}\varphi_i:\star_{i\in I}\mathcal{A}_{i}\rightarrow {\mathbb C}$ 
can be expressed in the form
$$
\star_{i\in I}\varphi_i (a_1\cdots a_n)=\sum_{w=j_1\cdots j_n\in \mathpzc{M}_n}\Phi(a_1(j_1)\cdots a_n(j_n))
$$
for any $a_1\in \mathcal{A}_{i_1}, \ldots, a_n\in \mathcal{A}_{i_n}$, where $i_1\neq \cdots \neq i_n$ and 
$a_k(j_k)$ is the orthogonal replica of $a_k$ of color $j_k$ for any $k$. 
The notation $\star_{i\in I}\mathcal{A}_{i}$ means that we take the free product of unital algebras 
with identified units in contrast to
$*_{i\in I}\mathcal{A}_{i}$ denoting the free product of unital algebras without identification of units.

It is noteworthy that a similar approach can be applied 
to describe nonsymmetric products of functionals, like the orthogonal and 
s-free products. This is a natural consequence of the fact that the orthogonal and s-free 
random variables can be expressed in terms of orthogonal replicas as shown in \cite{[17]}.  
A suitable family of path-dependent subordination functionals 
can be defined to decompose these products according to Motzkin paths.
The same holds for the monotone product of two functionals since in that case
the tensor product realization found in \cite{[L2019]} fits well into the replica scheme. 
We only give a glimpse of these results here and present them in more detail 
together with Motzkin decompositions of convolutions in a separate paper to keep the lenght of 
this one within reasonable limits. 
It is not clear whether orthogonal replicas can be used to obtain a Motzkin decomposition 
of the monotone product of an arbitrary totally ordered family of functionals, but 
we will also show there that slightly different replicas constructed in \cite{[L2019]} can certainly 
be used for that purpose. 

In standard approaches to noncommutative probability the information about distributions 
of variables is contained in a {\it sequence} of multivariate functionals, be it moments 
or cumulants (of one type or another), and in both cases 
the whole sequence of them needs to be studied together. 
The novelty in our approach is that we have a sequence of replicas of each variable, which, by a certain
duality relation, leads to a larger {\it family of path-dependent} multilinear functionals of rich combinatorial and algebraic 
structure which can also be studied together. On the combinatorial level, our approach fills the `gap' between 
interval partitions and noncrossing partitions of Kreweras \cite{[11]}. On the level of functionals, it
goes even deeper since it fills the `gap' between the boolean and free products of states on the algebraic level.

The main idea of this paper, to use Motzkin paths (or, words) as `arguments' of some basic functionals, has its roots 
in our construction of the {\it noncommutative logarithm of the Fourier transform} which 
was the ``boolean-classical logarithm of the Fourier transform'' based on the combinatorics on words \cite{[15]}.
This led us to the model in which ${\mathbb N}$ is replaced by $\mathpzc{M}$, equipped with its remarkable structure.
Our main tool is that of an algebraic decomposition of various functionals  
induced by the decomposition of free random variables into series of orthogonal replicas \cite{[14],[16]}
and by the decomposition of the free additive convolution \cite{[17]}. 
The main algebraic point of our approach on the level of variables is that the orthogonal replicas play the 
role of a `basis' of the `space of variables'. By a duality relation, the associated path-dependent functionals 
give a generating set of the space of moments and cumulants. This leads to additivity of functionals 
on the level of paths, involving moments and cumulants associated with different notions
of noncommutative independence. The latter and relations between them were studied by many authors \cite{[3],[M],[13],[22],[23],[24]}. The cumulants in our theory, called {\it Motzkin cumulants}, will be investigated in a forthcoming paper. 

The remaining part of the paper consists of 7 sections. In Section 2, we give some basic definitions and facts.
Orthogonal replicas and their moments are studied in Section 3. 
Lattices of Motzkin paths and Motzkin words are described in Section 4. 
In Section 5, we show that the moments of orthogonal replicas are described by Motzkin lattices, which leads to the concept of Motzkin 
functionals in terms of which the free and boolean products of functionals are 
decomposed. In Section 6, we prove important lemmas on the moments under Motzkin functionals.
In Section 7, we introduce and study the lattices $\mathcal{M}(w)$ of noncrossing partitions which are monotonically
adapted to $w\in \mathpzc{M}$, where monotonicity relates the depths of blocks to the heights of 
segments in Motzkin paths. Using $\mathcal{M}(w)$, we express in Section 8 the moments of Motzkin functionals in terms of boolean cumulants, which is a refinement of the formula for moments of free random variables in terms of boolean cumulants derived recently in \cite{[6], [10]}. 

\section{Preliminaries}

Moments and cumulants are functionals which contain information about distributions of random variables and 
the associated probability measures. In this paper, we concentrate on moments, although 
we formulate our results in a way that will enable us to treat the cumulants, using 
a similar language.

A noncommutative probability space is a pair $({\mathcal A},\varphi)$, where 
${\mathcal A}$ is a unital algebra and $\varphi$  is a normalized unital 
linear functional on ${\mathcal A}$ ($\varphi:{\mathcal A}\rightarrow {\mathbb C}$, $\varphi(1)=1$).
If ${\mathcal A}$ is a unital *-algebra, then we assume in addition that $\varphi$ is positive 
($\varphi(aa^*)\geq 0$ for any $a\in {\mathcal A}$). In both cases, we will sometimes call
such functionals {\it states}.
The {\it distribution} of $a\in \mathcal{A}$ with respect to $\varphi$ is the linear functional
$$
\mu_a:{\mathbb C}\langle X \rangle\rightarrow {\mathbb C}, \;\;\;
\mu_a(f)=\varphi(f(a)) \;\;{\rm for} \;\;f\in {\mathbb C}\langle a \rangle,
$$
where ${\mathbb C}\langle X \rangle$ is the unital algebra of polynomials in the indeterminate $X$.
Under certain assumptions, one can associate with this distribution a unique probability measure on the real line.
More generally, if $\{a_{i}\}_{i\in I}$ is a family of random variables from 
${\mathcal A}$, where $I$ is an index set, then the {\it joint distribution} 
of $\{a_{i}\}_{i\in I}$ with respect to $\varphi$ is the linear functional 
$$
\mu:{\mathbb C}\langle \{x_{i}\}_{i\in I}\rangle \rightarrow {\mathbb C}, \;\;\;
\mu(f)=\varphi(f(\{a_{i}\}_{i\in I}))\;\;{\rm for}\;\; f\in {\mathbb C}\langle \{x_{i}\}_{i\in I}\rangle ,
$$ 
where ${\mathbb C}\langle \{X_i\}_{i\in I}\rangle$ is the unital algebra of noncommutative polynomials
in the family of indeterminates $\{X_{i}\}_{i\in I}$.

\begin{Definition}
{\rm 
Let ${\mathcal A}$ be an algebra and let $\{f_n:\,n\in {\mathbb N}\}$ be a sequence of 
$n$-linear functionals $f_n:{\mathcal A}\times \cdots \times {\mathcal A}\rightarrow {\mathbb C}$
with values $f_n(a_1, \ldots , a_n)$. Introduce the associated family of $n$-linear functionals
\begin{eqnarray*}
f_{\pi_0}[a_1, \ldots , a_n]&:=&\prod_{V\in \pi_0}f(V)[a_1, \ldots , a_n],
\end{eqnarray*}
where $\pi_0\in P(n)$, $a_1, \ldots , a_n\in \mathcal{A}$, and
$f(V)[a_1, \ldots, a_n]:=f_k(a_{i(1)}, \ldots , a_{i(k)})$
for any set $V=\{i(1)< \cdots <i(k)\}\subset [n]$, where
$P(n)$ is the set of all partitions of $[n]$.
}
\end{Definition}

As we already mentioned in the Introduction, one of the main features of our approach is 
instead of indexing objects by natural numbers ${\mathbb N}$ we use the set of all nonempty reduced Motzkin words
$\mathpzc{M}$. In fact, we use a larger set of all (not necessarily reduced) Motzkin words $\mathpzc{AM}$ 
which begin and end with the same $h$, where $h\in \mathbb{N}$. For any 
$$
w=j_1\cdots j_n\in \mathpzc{AM}
$$ 
and any $\pi_{0}\in \mathcal{P}(n)$, where $n=|w|$, 
we assign the letter $j_p$ to the number $p\in [n]$ and thus, if 
$V_q=\{q(1), \ldots , q(k)\}\in \pi_0$, we assign the subword 
$$
v_q=j_{q(1)}\cdots j_{q(k)}
$$ 
to the block $V_q$. For simplicity, we will speak of the partition $\pi=\{v_1, \ldots , v_p\}$
associated with $\pi_0=\{V_1, \ldots, V_p\}$ and $w\in \mathpzc{AM}$, 
where we assume $V_1<\cdots <V_p$. It is natural to extend the above definition of multiplicative 
functionals to the context of partitions of Motzkin words $w$. For the sake of generality, 
in the definition given below we take an arbitrary set of words $\mathpzc{W}$ and, for any $w\in \mathpzc{W}$, 
the associated set of all partitions of $w$, denoted $\mathcal{P}(w)$. 
However, we will later use only Motzkin words and 
noncrossing partitions (adapted to $w$ in a suitable way).

\begin{Definition}
{\rm 
Let ${\mathcal A}$ be an algebra and let $\{f(w):\,w\in \mathpzc{W}\}$ be a family
of $|w|$-linear functionals $f(w):{\mathcal A}^{|w|}\rightarrow {\mathbb C}$
with values $f(w)(a_1, \ldots , a_n)$, where $\mathpzc{W}$ is a set of words.
Introduce the associated family of $n$-linear functionals
\begin{eqnarray*}
f_{\pi}[a_1, \ldots , a_n]&:=&\prod_{v\in \pi}f(v)[a_1, \ldots , a_n],
\end{eqnarray*}
where $\pi\in \mathcal{P}(w)$, $a_1, \ldots , a_n\in \mathcal{A}$, and 
$
f(v)[a_1, \ldots, a_n]:=f(v)(a_{i(1)}, \ldots , a_{i(k)})
$
for any subword $v=j_{i(1)}\cdots j_{i(k)}\in \mathpzc{W}$ of 
$w=j_1\cdots j_n\in \mathpzc{W}$.
}
\end{Definition}

Therefore, the family $\{f_{\pi}:\pi \in \mathcal{P}(w), w\in \mathpzc{W}\}$ can be treated as 
a multiplicative extension of the family $\{f(w):\,w\in \mathpzc{W}\}$, with 
multiplicativity induced by the set of subwords of $w$.
Note that the notation of Definition 2.1 is similar to that used in \cite{[22]} and \cite{[23]}, 
where parentheses are used for sequences of functionals and brackets indicate that we have functionals
related in one way or another to partitions. The only difference is that we write $f(w)$
instead of $f_{n}$. The main reason is that 
in Definition 2.2 we have words $w$ as arguments of $f(w)$ rather than $f_n$, where $n$ is merely an index, 
since these arguments will be subject to addition.

An important role in our computations is played by the {\it boolean cumulants}.
In this context, we use the lattices ${\rm Int}(n)$ of interval partitions of $[n]$ 
in the enumerative framework (see, for instance \cite{[24]}) and the associated lattices 
${\mathcal Int}(w)$ of interval partitions of $w\in \mathpzc{W}$ in the lattice framework.

\begin{Definition}
{\rm 
Let $\varphi$ be a linear normalized functional 
on the unital algebra $\mathcal{A}$.
By the {\it boolean cumulants} $\{\beta_n:\,n\in \mathbb{N}\}$ 
we understand the multilinear functionals defined by the recursive formula
$$
\varphi(z_1 \cdots z_n)=\sum_{\pi_{0}\in {\rm Int}(n)}\beta_{\pi_{0}}[z_1, \ldots, z_n],
$$
where 
$$
\beta_{\pi_{0}}[z_1, \ldots , z_n]=\prod_{V\in\pi_{0}}\beta_{|V|}[z_1, \ldots , z_n]
$$
for any $\pi_{0}\in {\rm Int}(n)$ and any $z_1, \ldots , z_n\in  \mathcal{A}$. 
}
\end{Definition}

In this connection, let us observe that there is a large `combinatorial gap' between boolean cumulants and free 
cumulants involving noncrossing partitions in the sense that for large $n$ the set ${\rm Int}(n)$ is much smaller than 
${\rm NC}(n)$ and thus the underlying combinatorics of moments and cumulants 
is much simpler. One of our original motivations was to investigate this `combinatorial gap' on an algebraic level.
We find it interesting that it is filled with the structure of Motzkin paths which 
show up on the level of moments and cumulants.  

Let us remark that we will usually deal with mixed moments, or even mixed moment functionals, 
but we will call them moments for simplicity if no confusion arises.

\section{Orthogonal replicas}

Let us present some details on the decomposition of free random variable in terms of orthogonal replicas.
For more details, see \cite{[14],[16],[17]}.

For simplicity, consider first two noncommutative probability spaces, $({\mathcal A}_{i}, \varphi_i)$, 
where $i\in J=\{1,2\}$.
\begin{enumerate}[(a)]
\item
Let $p$ be an abstract projection, namely $p^2=p$ (in the case of the *-algebra we also 
assume that $p=p^*$). For any noncommutative probability space
$({\mathcal A}, \varphi)$, we first construct its $p$-{\it extension} 
$(\widetilde{\mathcal A}, \widetilde{\varphi})$, where 
$$
\widetilde{\mathcal A}={\mathcal A}\star{\mathbb C}[p]
$$ 
is the free product with identified units (with the natural involution if needed) and
$\widetilde{\varphi}$ is the linear extension of 
$$
\widetilde{\varphi}(p^{\gamma}a_1pa_2p\cdots pa_mp^{\delta})
=
\varphi(a_1)\varphi(a_2)\cdots \varphi(a_m)
$$  
and $\widetilde{\varphi}(p)=1$, where $\gamma, \delta\in \{0,1\}$ and 
$a_1,\ldots , a_m\in {\mathcal A}$. In the case of a *-algebra, if $\varphi$ is a state, 
then $\widetilde{\varphi}$ is also a state. This is a boolean way of extending $\varphi$ \cite{[14]}. Clearly, the distribution (*-distribution) of any $a\in {\mathcal A}$ in 
the states $\widetilde{\varphi}$ and $\varphi$ are identical.
\item
Let $\widetilde{\mathcal{A}}$ be a $p$-extension of ${\mathcal A}$.
Construct the infinite tensor product $\widetilde{\mathcal A}^{\otimes \infty}$ as follows. 
First, form the linear vector space ${\mathcal F}$ 
with the basis given by elements $(x_j)_{j\in \mathbb{N}}$ of the Cartesian product 
$\prod_{j\in {\mathbb N}}\widetilde{{\mathcal A}}$ with the property that almost all 
coordinates are in the set $G:=\{1,p\}$.
Then ${\mathcal F}$ becomes an algebra if we extend linearly the
multiplication for basis elements given by
$$
(x_{j})_{j\in \mathbb{N}}(x_{j}')_{j\in \mathbb{N}}=(x_{j}x_{j}')_{j\in \mathbb{N}}
$$
to the vector space ${\mathcal F}$. Let now ${\mathcal F}_{0}$ denote the vector subspace
of ${\mathcal F}$ spanned by its elements of the form
$$
(x_{j})_{j\in \mathbb{N}}-(x'_{j})_{j\in \mathbb{N}}-\alpha (x_{j}'')_{j\in \mathbb{N}},
$$
where
$$
x_{j}=x_{j}'=x_{j}'' \;\;{\rm for}\;\; j\neq j_{0}\;\; {\rm and}\;\;
x_{j_{0}}=x_{j_{0}}'+\alpha x_{j_{0}}''
$$
where $j_{0}$ runs through the index set ${\mathbb N}$ and 
$x_{j},x_{j_0},x_{j_0}', x_{j_0}''$ run through the algebra $\widetilde{\mathcal A}$
and $\alpha \in {\bf C}$. Since ${\mathcal F}_{0}$ is an ideal in ${\mathcal F}$, we can form
the  quotient algebra
$$
\bigotimes_{j\in {\mathbb N}}\widetilde{\mathcal A}\equiv \widetilde{\mathcal A}^{\otimes \infty}={\mathcal F}/{\mathcal F}_{0}
$$
which is the infinite tensor power of the
algebra $\widetilde{\mathcal{A}}$ with respect to the set $G$.
The elements of the quotient basis are obtained in a canonical fashion
from the basis of $\mathcal{F}$ and are denoted by $\bigotimes_{j\in {\mathbb N}}x_{j}$.
This is a special case of the infinite tensor product of algebras with respect to a
family of elements given in \cite{[FLS]}.
\item
Let $\widetilde{\mathcal{A}}_1$ be the $p$-extension of ${\mathcal A}_1$,
and let $\widetilde{\mathcal{A}}_2$ be the $q$-extension of ${\mathcal A}_2$.
In the tensor product
$$
{\mathcal A}_{\otimes}:=
\widetilde{\mathcal A}_{1}^{\otimes \infty}\otimes \widetilde{\mathcal A}_{2}^{\otimes\infty}
$$
we distinguish families of orthogonal projections, 
$\{p(j):j\in \mathbb{N}\}$ and $\{q(j):s\in \mathbb{N}\}$, defined as follows: $p(1)=p^{\otimes \infty}$, 
$q(1)=q^{\otimes \infty}$,  
and
\begin{eqnarray*}
p(j)&=&1_{1}^{\otimes (j-2)}\otimes p^{\perp}\otimes p^{\infty},\\
q(j)&=&1_{2}^{\otimes (j-2)}\otimes q^{\perp}\otimes q^{\infty},
\end{eqnarray*}
for $j>1$, where $1_i$ is the unit in ${\mathcal A}_{i}$. It is convenient to set $q=p$ 
in proofs and examples. 
\item
For any $a\in {\mathcal A}_{1}$ and $b\in {\mathcal A}_{2}$, we shall use canonical 
injections into $\widetilde{\mathcal{A}}_{1}^{\otimes \infty}$
and $\widetilde{\mathcal A}_{2}^{\otimes \infty}$, respectively, namely
\begin{eqnarray*}
\iota_j(a)&=&1_{1}^{\otimes (j-1)}\otimes a\otimes 1_{1}^{\otimes \infty},\\
\iota_j(b)&=&1_{2}^{\otimes (j-1)}\otimes b\otimes 1_{2}^{\otimes \infty},
\end{eqnarray*}
for any natural $j$.
\end{enumerate}

\begin{Definition}
{\rm Let $(\widetilde{\mathcal A}_{1}, \widetilde{\varphi}_{1})$ be a $p$-extension of 
$({\mathcal A}_1, \varphi_1)$ and let $(\widetilde{\mathcal A}_{2}, \widetilde{\varphi}_{2})$ 
be a $q$-extension of $({\mathcal A}_2, \varphi_2)$. Consider the noncommutative probability space 
$({\mathcal A}_{\otimes}, \Phi_{\otimes})$, where 
$$
\Phi_{\otimes}:=\Phi_1\otimes \Phi_2\;\;\;{\rm and}\;\;\; {\mathcal A}_{\otimes}:=
\widetilde{\mathcal A}_{1}^{\otimes \infty}\otimes \widetilde{\mathcal A}_{2}^{\otimes\infty}
$$
and $\Phi_j=\widetilde{\varphi}_{j}^{\otimes\infty}$ for $j=1,2$. The elements
\begin{eqnarray*}
a(j)&:=&\iota_{j}(a)\,\otimes \,q(j),\\
b(j)&:=&p(j)\,\otimes\, \iota_j(b),
\end{eqnarray*}
where $a,a'\in \mathcal{A}_{1}, b,b'\in \mathcal{A}_2$ and $j\in {\mathbb N}$, are called {\it orthogonal replicas} (or, {\it replicas})
of $a$ and $b$, respectively. The noncommutative probability space $({\mathcal A}_{{\rm rep}}, \Phi)$, where
${\mathcal A}_{{\rm rep}}$ is the unital subalgebra of ${\mathcal A}_{\otimes}$ generated by orthogonal replicas 
and $\Phi$ is the restriction of $\Phi_{\otimes}$ to ${\mathcal A}_{{\rm rep}}$, will be called the {\it orthogonal replica space}
(or, the {\it replica space}).}
\end{Definition}

\begin{Remark}
{\rm 
Let us make some comments on the above definition.
\begin{enumerate}[(a)]
\item
We call these replicas orthogonal since
$$
a(j)a'(k)=b(j)b'(k)=0
$$
whenever $j\neq k$, for any $a\in \mathcal{A}_1$ and $b\in \mathcal{A}_2$.
\item
The orthogonal replicas of units $1_1\in {\mathcal A}_{1}$ and $1_2\in {\mathcal A}_{2}$ do not coincide with the unit
in the orthogonal replica space, which is important in any unification 
context since, for instance, boolean independence does not admit unification of units. 
For that reason, subalgebras of the orthogonal replica space 
which are generated by some (but not all) orthogonal replicas of some elements from 
${\mathcal A}_{1}$ and some elements from ${\mathcal A}_{2}$ are not unital, although they may have
internal units.
However, one can construct the unit in the orthogonal replica space ${\mathcal A}_{{\rm rep}}$ from these replicas since  
$$
\sum_{j=1}^{\infty}1_1(j)=\sum_{j=1}^{\infty}1_2(j)=1_{{\rm rep}} \;,
$$
where $1_{{\rm rep}}$ is the unit in ${\mathcal A}_{{\rm rep}}$, which gives the unit identification
when needed (like in the free case).
\item
If $({\mathcal A}_{1}, \varphi_1)$ and $({\mathcal A}_{2}, \varphi_2)$
are *-probability spaces or $C^*$-probability spaces and we perform the GNS construction for them, 
we obtain triples $({\mathcal H}_{1}, \pi_1, \xi_1)$ and $({\mathcal H}_{2}, \pi_2, \xi_2)$, 
where ${\mathcal H}_{i}$ is a pre-Hilbert or Hilbert space, respectively, $\pi_i$ is a *-representation 
and $\xi_i$ is a cyclic unit vector. Then, it is not hard to see that the GNS triples
for $(\widetilde{\mathcal A}_{1}, \widetilde{\varphi}_{1})$ and $(\widetilde{\mathcal A}_{2}, \widetilde{\varphi}_{2})$ are
$({\mathcal H}_{1}, \widetilde{\pi}_{1}, \xi_{1})$ and $({\mathcal H}_{2}, \widetilde{\pi}_{2}, \xi_{2})$,
where $\widetilde{\pi}_{1}$ and  $\widetilde{\pi}_{2}$ are the unique extensions of $\pi_1$ and $\pi_2$ 
which map $p$ and $q$ onto the projections onto ${\mathbb C}\xi_1$ or ${\mathbb C}\xi_2$, respectively. 
\item
Let us remark that tensor product realizatons show how 
to associate products of graphs with various concepts. For instance, if $a$ and $b$ are adjacency matrices 
of rooted graphs $({\mathcal G}_{1},e_1)$, $({\mathcal G}_{2}, e_2)$, the variables $a(1)$ and $b(1)$ are 
boolean independent under $\Phi$ and their moments reduce to those of
$$
a\otimes p\;\;\;{\rm and}\;\;\; q\otimes b.
$$
Then we can interpret this realization as follows: 
$p$ indicates that we should glue ${\mathcal G}_{1}$ to the root $e_2$ of ${\mathcal G}_1$ and 
$q$ indicates that we should glue ${\mathcal G}_{2}$ to the root $e_1$ of ${\mathcal G}_2$.
Other graph products can be treated in a similar fashion \cite{[1]}.

\end{enumerate}
}
\end{Remark}

Below, we shall present a useful application of the orthogonal projections, which gives a concise formula for boolean cumulants. 
Moments of the form given below were studied in Proposition 5.1 in \cite{[17]}, where it was shown that they were boolean 
cumulants of one variable (more precisely, we used there $(-1)^{n}\beta_n$, called `inverse boolean cumulants'). The multivariate case is analogous.

\begin{Proposition}
Let $a_1, \ldots , a_n$ be elements of a noncommutative probability space 
$({\mathcal A}, \varphi)$ and let $\widetilde{\varphi}$ be the extension of $\varphi$ to ${\mathcal A}(p)$.
Then, it holds that
$$
\widetilde{\varphi}(a_{1}p^{\perp}a_{2}\cdots p^{\perp}a_{n})
=\beta_{n}(a_1,a_2, \ldots ,a_n),
$$
where $\beta_{n}(a_1, a_2, \ldots , a_n)$ is the mixed boolean cumulant of $a_1, a_2, \ldots , a_n$ of order $n$.
\end{Proposition}
{\it Proof.}
The proof is elementary and follows from that of Proposition 5.1 in \cite{[17]}.
\hfill $\blacksquare$\\

\begin{Example}
{\rm Let us compute a simple boolean cumulant, using Proposition 3.1:
\begin{eqnarray*}
\widetilde{\varphi}(a_1p^{\perp}a_2p^{\perp}a_3)&=& \widetilde{\varphi}(a_1a_2a_3)-\widetilde{\varphi}(a_1pa_2a_3)-
\widetilde{\varphi}(a_1a_2pa_3)+\widetilde{\varphi}(a_1pa_2pa_3)\\
&=&{\varphi}(a_1a_2a_3)-\varphi(a_1)\varphi(a_2a_3)-\varphi(a_1a_2)\varphi(a_3)+\varphi(a_1)\varphi(a_2)\varphi(a_3)\\
&=& \beta_3(a_1,a_2,a_3)
\end{eqnarray*}
}
\end{Example}

\begin{Corollary}
If $[n]=J_1\cup \ldots \cup J_p$ is a partition into nonempty 
disjoint intervals and $b_j=\prod_{k\in J_j}a_k$ for any $j=1, \ldots , r$ (with increasing order 
of indices), then 
$$
\widetilde{\varphi}(b_1p^{\perp}b_{2}\cdots p^{\perp}b_{n})
=\sum_{\stackrel{\pi\in {\rm Int}(n)}
{\scriptscriptstyle J_1\sim \cdots \sim J_{r}}}
\beta_{\pi}[a_1,\ldots ,a_n],
$$
where notation $J_j\sim J_{j+1}$ means that $J_j$ is connected with $J_{j+1}$ by a block of $\pi$.
\end{Corollary}
{\it Proof.}
This is a straightforward consequence of Proposition 3.1 and the moment-cumulant formula for boolean cumulants.
\hfill $\blacksquare$

\begin{Example}
{\rm Let us compute a simple example:
\begin{eqnarray*}
\widetilde{\varphi}(a_1p^{\perp}a_2a_3)&=&
\beta_2(a_1,a_2a_3)= 
\varphi(a_1a_2a_3)-\varphi(a_1)\varphi(a_2a_3)\\
&=&
\beta_3(a_1,a_2,a_3)+\beta_2(a_1,a_2)\beta_1(a_3)+ \beta_1(a_1)\beta_2(a_2,a_3)\\
&+& \beta_1(a_1)\beta_1(a_2)\beta_1(a_3)-\beta_1(a_1)(\beta_2(a_2,a_3)+\beta_1(a_2)\beta_1(a_3))\\
&=&
\beta_3(a_1,a_2,a_3)+\beta_2(a_1,a_2)\beta_1(a_3)
\end{eqnarray*}
}
\end{Example}

The main motivation to study orthogonal replicas is that they show up
in the decomposition of free random variables studied in our previous papers.
The case of two algebras is especially illuminating, for which we state the main result of 
\cite{[14]} and \cite{[16]}, using the notation introduced in this paper.

\begin{Theorem}
Let $({\mathcal A}_{1}, \varphi_1)$ and $({\mathcal A}_{2}, \varphi_{2})$ 
be noncommutative probability spaces and let $K_1\subset {\mathcal A}_{1}$ and $K_2\subset {\mathcal A}_{2}$.
Then the families of the form
$$
\{\sum_{j=1}^{\infty}a(j): a\in K_1\}\;\;\;and\;\;\; 
\{\sum_{j=1}^{\infty}b(j): b\in K_2\}
$$
are free with respect to $\Phi$. Moreover, they have 
have the same joint distributions with respect to $\Phi$ as the
families $K_1$ and $K_2$ with respect to $\varphi_1$ and $\varphi_2$, respectively.
\end{Theorem}

This result was proved in \cite{[14],[16]} in two versions: first, sequences of random variables which approximate free random variables were constructed in \cite{[14]} and then a meaning to infinite series realization was given in \cite{[16]}. 
Below we state that result, although we refrain from using 
the symbol $\overline{\otimes}$ introduced in \cite{[16]} to implement
Berberian's theory \cite{[5]} to deal with infinite series. For the purposes of this paper, it 
suffices to think of the series as sequences of partial sums since for each moment only a finite number of summands suffices 
to reproduce the moments of free random variables, the remaining summands give zero contribution. On the dual level of functionals,
and that approach we pursue, it is even more apparent since each of the functionals is well-defined pointwise.

Orthogonal, s-free and monotone random variables can be treated in a similar fashion in the case of two 
noncommutative probability spaces. Namely, they can be expressed in terms of 
orthogonal replicas. This follows from \cite{[17]} for the orthogonal and 
s-free random variables and from \cite{[L2019]} for monotone random variables. 
Let us remark that the tensor product 
realization of monotone random variables in the spirit of \cite{[14]} given by Franz \cite{[7]} 
is not done in terms of orthogonal replicas and therefore it is probably not appropriate 
for the approach via Motzkin paths.

\begin{Theorem}
Let $({\mathcal A}_{1}, \varphi_1),({\mathcal A}_{2}, \varphi_2)$ be 
noncommutative probability spaces and let 
$K_1\subset {\mathcal A}_{1}$ and $K_2\subset {\mathcal A}_{2}$. Then
\begin{enumerate}
\item
the families of the form
$$
\{\sum_{j\;\;{\rm odd}}a(j):a\in K_{1}\}\;\;\;and\;\;\;\{\sum_{j\;\;{\rm even}}b(j):b\in K_{2}\}
$$
are s-free with respect to $\Phi$,
\item
the families of the form
$$
\{a(1):a\in K_1\}\;\;\;and \;\;\;\{b(2): b\in K_2\},
$$
are orthogonal with respect to $\Phi$,
\item
the families of the form
$$
\{a(1): a\in K_1\}\;\;\;and\;\;\;\{b(1)+b(2): b\in K_2\}
$$
are monotone independent with respect to $\Phi$.
\end{enumerate}
In all these cases, the given families have the same joint distributions
with respect to $\Phi$ as the families $K_1$ and $K_2$ with respect to $\varphi_1$ and $\varphi_2$, respectively.
\end{Theorem}
{\it Sketch of the proof.}
The Hilbert space realizations of orthogonal and s-free random variables 
in terms of orthogonal replicas were given in Theorem 4.2 and 
Theorem 7.1 in \cite{[17]} for $C^*$-probability spaces. 
A realization of monotone (in fact, even conditionally monotone) random variables 
in terms of orthogonal replicas for noncommutative probability 
spaces was given in \cite{[L2019]}. Note that in the case of two noncommutative probability 
spaces the realizations given there reduce to orthogonal replicas, but in the general case 
these realizations are in terms of slightly different replicas. We treat these versions of 
proofs are sufficient for our present purposes. Detailed algebraic proofs will be given in a separate paper.
\hfill $\blacksquare$\\

\begin{Remark}
{\rm The tensor realization of Theorem 3.1 can be generalized to the case of 
an arbitrary family of noncommutative probability spaces $\{({\mathcal A}_{i}, \varphi_{i}):i\in I\}$, where 
$I$ is an arbitrary index set, 
except that one needs to introduce slightly more general families of orthogonal projections.
An outline of the general case is given below.
\begin{enumerate}[(a)]
\item
Let  $\{p_{i}(j):j\in \mathbb{N}\}$ be a sequence of orthogonal projections for each $i\in I$
which act identically onto tensor sites associated with labels different form $i$, namely 
$p_i(1)=\bigotimes_{k\neq i}p_{k}^{\otimes \infty}$ for any $i$ and 
$$
p_i(j)=P_i(j)-P_i(j-1), \;\;\;
{\rm where} \;\;\;P_i(j)=\bigotimes_{k\neq i}(1_{k}^{\otimes (j-1)}\otimes p_k^{\otimes \infty}),
$$
for $n>1$, where $p_k=p_{k}^2$ is used to construct the $p_k$-extension of $\varphi_{k}$, and we set
$P_i(0)=0$ for any $i$. 
\item
For any $a\in {\mathcal A}_{k}$, let 
$$
a(j)=\iota_{j}(a)\otimes p_k(j),
$$
where $\iota_j(a)$ is a copy of $a$ of the form
$$
\iota_j(a)=1_{k}^{\otimes (j-1)} \otimes a\otimes 1_{k}^{\otimes \infty}.
$$
We will call $a(j)$ an {\it orthogonal replica} (or, {\it replica}) of $a$, with 
{\it label} $k$ and {\it color} $j$, respectively.
If we vary $j$, the positions of variables with 
label $k$ in the tensor product change accordingly. It was only for the sake of convenience that
our variable with label $k$ was written above as the first one in the tensor product.
\item
Observe that if some variables from ${\mathcal A}_{i}$ satisfy a relation, 
then an analogous relation holds for their replicas of label $i$ and 
color $j$ for fixed $i$ and $j$. In fact, the mapping $\tau:{\mathcal A}_{i}\rightarrow {\mathcal A}_{{\rm rep}}$ 
given by 
$$
\tau(a)=a(j)
$$ 
is a non-unital homomorphism for any $i$ and $j$. Of course, the unit in the replica space is of the form
$$
1_{{\rm rep}}=\bigotimes_{i\in I}1_{i}^{\otimes \infty}
$$ 
and if $1_i$ is a unit in ${\mathcal A}_{i}$, then each 
$$
1_{i}(j)=1_{i}^{\otimes \infty}\otimes p_i(j)\neq 1_{{\rm rep}}
$$ 
for any $j$. However, each family $\{1_{i}(j):\,j\in \mathbb{N}\}$ gives a decomposition of $1_{{\rm rep}}$.
\item
The noncommutative probability space of interest, in which the above variables live, 
is of the form $({\mathcal A}_{{\rm rep}}, \Phi)$, where
$$
{\mathcal A}_{\otimes}=\bigotimes_{i\in I}\widetilde{{\mathcal A}_{i}}^{\otimes \infty}\;\;\;
{\rm and}\;\;\;
\Phi_{\otimes}=\bigotimes_{i\in I}\widetilde{\varphi}_{i}^{\otimes \infty},
$$
where we assume that ${\mathcal A}$ is the linear span of tensors
$\bigotimes_{i\in I}\bigotimes_{k=1}^{\infty} a_{i,k}$ such that $a_{i,k}\notin\{1_{i},p_{i}\}$ for 
a finite number of pairs $(i,k)$.
By $\mathcal{A}_{{\rm rep}}$ we denote the unital algebra generated by all replicas
and $\Phi$ will stand for the restriction of $\Phi_{\otimes}$ to $\mathcal{A}_{{\rm rep}}$.
Details of this approach for an arbitrary index set $J$ can be found in \cite{[14],[16]}.
\end{enumerate}
}
\end{Remark}

We would like to study mixed moments of orthogonal replicas in more detail and show that there is 
a Motzkin path structure behind their combinatorics that we have discovered recently.  
In order to state the main results involving the Motzkin structure, let us look closer at
the moments of replicas. Since the latter give decompositions of free random variables, it is not surprising that if we take replicas of variables which are in the kernels of the corresponding states and neighbors belong to different algebras, then their mixed moments vanish. Thus, we obtain a kind of freeness condition. However, one has to remember that the replicas of the associated units are not identified. Therefore, one has to provide additional recursions for moments involving the replicas of units.
Let us first state some basic orthogonality results.
\begin{Proposition}
Let  $a,a'\in {\mathcal A}_{i}$, $b\in {\mathcal A}_{l}$, $a_{k}\in \mathcal{A}_{i_{k}}$ where
$i\neq i_1\neq \cdots \neq i_n\neq l$, and let $z(w):=a_{1}(j_1) \cdots a_{n}(j_n)$,
where  $w=j_1\cdots j_n$ and $j_k>j$ for all $1\leq k\leq n$, $n\in \mathbb{N}$. We have:
\begin{enumerate}
\item
orthogonality: if $j\neq j'$, then 
$$
a(j)a'(j')=0,
$$ 
\item
remote orthogonality: if $i\neq l$ and $j_k=j+1$ for some $k$, then
$$
a(j)z(w)b(j)=0,
$$
\item
conditional orthogonality: if $i=l=i_k$ and $j_k=j+1$ for some $k$, then
$$
a(j)z(w)b(j)=0.
$$
\end{enumerate}
\end{Proposition}
{\it Proof.}
The orthogonality of replicas which have the same labels 
follows from the fact that if $j\neq k$ then 
$p_{i}(k)\perp p_{i}(l)$ for $k\neq l$ and arbitrary $i\in I$.
To prove the remote orthogonality of replicas which have different labels, one has to 
find a place where some cyclic projection $p_q$ meets $p_{q}^{\perp}$.
Let us look at sites of color $j$. 
There is at least one replica in $z(w)$ with a projection of type $p_{r}(j+1)$ which can be decomposed 
into a sum of a finite number of mutually orthogonal projections, each of which has either $p_i^{\perp}$ at site $(i,j)$ 
or $p_l^{\perp}$ at site $(l,j)$ and, moreover, the first one meets $p_l$ on the left, produced by $a(j)$, 
or the second one meets $p_i$ on the right, produced by $b(j)$. 
This means that all these terms vanish. To show more explicitly how this argument works, 
let us consider the product of three replicas $a(1)z(2)b(1)$. The essential part of this product (we restrict 
our attention to active tensor sites) is equal to
$$
\big((a\otimes 1_i)\otimes (p_l\otimes p_l)\otimes (p_r\otimes p_r)\big)
z(2)
\big((p_i\otimes p_i)\otimes (b\otimes 1_l)\otimes (p_r\otimes p_r)\big),
$$
where
$$
z(2)=\big((1_i\otimes p_i)\otimes (1_l\otimes p_l)-
(p_i\otimes p_i)\otimes (p_l\otimes p_l)\big)\otimes (1_r\otimes z).
$$
Decomposing $1_i=p_i+p_i^{\perp}$ and $1_l=p_l+p_l^{\perp}$ in the last expression, 
we can see that $z(2)$ 
is a tensor product of $(1_r\otimes z)$ with a 
sum of three orthogonal projections:
$$
(p_{i}^{\perp}\otimes p_i)\otimes (p_{l}^{\perp}\otimes p_l),\;\;
(p_{i}^{\perp}\otimes p_i)\otimes (p_{l}\otimes p_l),\;\;
(p_{i}\otimes p_i)\otimes (p_{l}^{\perp}\otimes p_l)
$$
and in each of these projections there is either $p_i^{\perp}$ that 
meets $p_i$ on the right or $p_l^{\perp}$ that meets $p_l$ on the left, which makes the product vanish.
The proof in the general case is very similar.
This proves (2). The proof of (3) is very similar to that of (2)
and is left to the reader. This completes the proof of the lemma.
\hfill $\blacksquare$

\begin{Lemma}
Let $a_k\in {\mathcal A}_{i_k}$, where $k=1, \ldots , n$
and $i_1\neq \cdots \neq i_n$. Then we have 
\begin{enumerate}
\item
the freeness property:
$$
\Phi\left(a_1(j_1)\cdots a_{n}(j_n)\right)=0
$$
whenever $a_k\in {\rm Ker}(\varphi_{i_k})$ for $k=1, \ldots n$, and 
$j_1, \ldots , j_n\in \mathbb{N}$,
\item
the reduction property:
$$
\Phi\left(a_1(j_1)\cdots 1_{i_r}(j_r) \cdots a_{n}(j_n)\right)
=
\Phi\left(a_1(j_1)\cdots \check{1}_{i_r}(j_r) \cdots a_{n}(j_n)\right)
$$
whenever $a_k\in {\rm Ker}(\varphi_{i_k})$ for
$k=1, \ldots r-1$, and $(j_1, \ldots , j_r)=(1,\ldots , r)$, 
where the element with $\check{}$ on top has to be omitted, 
and the moment on the LHS vanishes for the remaining colors.
\end{enumerate}
\end{Lemma}
{\it Proof.}
Let us adopt the convention for tensor products $\bigotimes_{i\in I}\bigotimes_{k=1}^{\infty}a_{i,k}$ 
that we write explicitly only $a_{i,k}\notin \{1_i,p_i\}$. 
Moreover, if we have such $a_{i,k}=a_q\in \mathcal{A}_i$ for some element $a_q\notin \{1_i,p_i\}$, we denote 
it $(a_q)_{i,k}$ to indicate the label and the color of $a_q$. This notation allows us to write
elements of the replica space in any order.
Therefore, a typical element of ${\mathcal A}_{{\rm rep}}$ will be of the form
$$
(a_{1})_{i_1,j_1}\otimes \cdots \otimes (a_n)_{i_n,j_n}
$$
where $a_k\in \mathcal{A}_{i_k}$, $k=1, \ldots , n$. We will also denote
$$
P=\bigotimes_{i\in I}p_i^{\infty}.
$$
It is convenient to treat $p_i$ as a projection onto the cyclic subspace ${\mathbb C}\xi_i$ even 
if we work in the purely algebraic setting. Let us assume now that $a_k\in \mathcal{A}_{i_k}\cap {\rm Ker}(\varphi_{i_k})$ 
for $k=1, \ldots n$. Now, it is clear that 
$$
Pa_1(j_1)\cdots a_n(j_n)P\equiv a_1(j_1)\cdots a_n(j_n)
$$ 
for any $j_1, \ldots, j_n$, where we write $x\equiv y$ if $\Phi(x)=\Phi(y)$. Next,
\begin{eqnarray*}
Pa_1(j_1)&= &0 \;\;\;{\rm unless}\;\;\; j_1=1\\
Pa_1(j_1)a_2(j_2)&= &0 \;\;\;{\rm unless}\;\;\;j_1=1,j_2=2\\
&\ldots&\\
Pa_1(j_1)\cdots a_{n}(j_n)&= &0 \;\;\;{\rm unless}\;\;\;j_1=1, j_2=2, \ldots , j_n=n
\end{eqnarray*}
since $P1_{i}(1)=P$ and $P1_{i}(j)=0$ for any $j>1$ and any $i$, and since
$Pa_{1}(1)a_{2}(j_2)=0$ if $j_{2}\neq 2$, $Pa_1(1)a_2(2)a_{3}(j_3)= 0$ 
if $j_3\neq 3$, etc. These equations follow from the fact that $p_i(j)$ for $j>1$ project 
onto the orthogonal complements of cyclic subspaces for labels $k\neq j$ and colors equal to $j-1$.
This gives the consecutiveness of colors as a necessary condition for the moment not to vanish.
A similar property holds if we multiply our product of variables by $P$ from the right.
In particular,
$$
a_{n}(j_n)P\equiv 0\;\;\;{\rm unless}\;\;\;j_n=1,
$$
and thus 
$$
Pa_1(j_1)\cdots a_{n}(j_n)P\equiv 0 
$$ 
for $n>1$. In turn, if $n=1$, we also obtain zero since $a_1\in {\rm Ker}(\varphi_{i_1})$. 
This proves the freeness property. 
In order to prove the reduction property, it suffices to consider
$$
Pa_1(j_1)\cdots a_{r-1}(j_{r-1})\equiv (a_1)_{i_1,j_1}\otimes \cdots \otimes (a_{r-1})_{i_{r-1},j_{r-1}}
$$
for $(j_1, \ldots, j_{r-1})=(1,\ldots, r-1)$ and multiply it from the right by $1_{i_r}(j_r)$ since, 
by the previous argument, this product vanishes for the remaining color sequences. Now,
$$
\left((a_1)_{i_1,1}\otimes \cdots \otimes (a_{r-1})_{i_{r-1},r-1}\right) 1_{i_r}(j_r)\equiv 0
$$
when $j_r\leq r-1$ since $p_{i_r}(j_r)$ contains the projection $p_{i_{r-1}}$ 
which meets $a_{r-1}\in {\rm Ker}(\varphi_{i_{r-1}})$ at site 
of color $r-1$. On the other hand, if $j_r> r$, then the above product is also 
equal to zero since in this case all sites of colors $\geq r$ are occupied by units or cyclic 
projections and $p_{i_r}(j_r)$ is a difference of two projections. Finally,
$$
a_{r-1}(r-1)1_{i_r}(r)=a_{r-1}(r-1)
$$
since $a_{r-1}\in {\rm Ker}(\varphi_{i_{r-1}})$. This shows the reduction property and thus the 
proof is completed.
\hfill $\blacksquare$

\begin{Corollary}
Let $a_k\in {\mathcal A}_{i_k}$, where $k=1, \ldots , n$
and $i_1\neq \cdots \neq i_n$. Then we have
\begin{enumerate}
\item the multiple reduction property:
$$
\Phi\left(1_{i_1}(j_1)\cdots 1_{i_r}(j_r) a_{r+1}(j_{r+1})\cdots a_{n}(j_n)\right)
=
\Phi\left(a_{r+1}(j_{r+1})\cdots a_{n}(j_n)\right)
$$
whenever $j_1=\cdots =j_r=1$ and $r\geq 1$, and otherwise this moment vanishes.
\item
the reversed reduction property:
$$
\Phi\left(a_1(j_1)\cdots 1_{i_r}(j_r) \cdots a_{n}(j_n)\right)
=
\Phi\left(a_1(j_1)\cdots \check{1}_{i_r}(j_r) \cdots a_{n}(j_n)\right)
$$
whenever $a_k\in {\rm Ker}(\varphi_{i_k})$ 
$k=r+1, \ldots n$, and $(j_n, \ldots, j_{r+1})=(1, \ldots, n-r+1)$, and for other values of the labels
the moment on the LHS vanishes.
\end{enumerate}
\end{Corollary}
{\it Proof.}
The multiple reduction property follows from a repeated application of the reduction property (for $r=1$).
The reversed reduction property is a mirror reflection of the reduction property and is proved in the same way.
\hfill $\blacksquare$\\

\section{Lattices $\mathpzc{M}_{n}$}

An important feature of the tensor product realization of free random variables given in 
\cite{[14],[16]} is that the computations of mixed moments of free random variables reduce
to the computations of mixed moments of copies of these variables whose colors 
form reduced Motzkin words when multiplied. This will be stated and shown in this work since 
at the time of writing the above-mentioned papers we were not aware of this fact. However, we used words to encode 
information on mixed moments in a related model that connected boolean and classical cumulants \cite{[15]}
and the present paper is a continuation of this approach.

\begin{Definition}
{\rm 
By a {\it reduced Motzkin word} of lenght $n\in \mathbb{N}$ we will understand a word in letters from the alphabet ${\mathbb N}$
of the form
$$
w=j_{1}j_{2}\cdots j_{n},
$$
where $j_{1}=j_{n}=1$ and $|j_{i}-j_{i-1}|\in \{-1,0,1\}$. 
The set of reduced Motzkin words of lenght $n$ will be denoted by $\mathpzc{M}_{n}$. We denote by 
$$
\mathpzc{M}=\bigcup_{n=1}^{\infty}\mathpzc{M}_{n}
$$ 
the set of all nonempty reduced Motzkin words. By a {\it Motzkin word} of lenght $n\in \mathbb{N}$ we shall understand 
a word of the form 
$$
w=j_1j_2\cdots j_n,
$$
where $j_1=j_n=h$ for some $h\in {\mathbb N}$ and $j_k\geq j_1$ for all $k$, with
$|j_{i}-j_{i-1}|\in \{0,1\}$. 
We will write in this case $h(w)=h$ and say that the {\it height} of $w$ is $h$. 
In particular, $h(j)=j$. 
The set of all Motzkin words will be denoted $\mathpzc{AM}$. We have
$$
\mathpzc{AM}=\bigcup_{n=1}^{\infty}\mathpzc{AM}_{n},
$$ 
where $\mathpzc{AM}_{n}$ is the lattice of all Motzkin words of lenght $n$. 
We denote $\mathpzc{M}^{*}=\mathpzc{M}\cup\{\varnothing\}$ 
and $\mathpzc{AM}^{*}=\mathpzc{AM}\cup\{\varnothing\}$, where $\varnothing$ is the empty word.
}
\end{Definition}

\begin{Example}
{\rm 
We list all reduced Motzkin words of lenght smaller than six:

\begin{eqnarray*}
\mathpzc{M}_{1}&=&\{1\},\\
\mathpzc{M}_{2}&=&\{11\},\\
\mathpzc{M}_{3}&=&\{111, 121\}, \\
\mathpzc{M}_{4}&=&\{1111, 1121, 1211, 1221\} \\
\mathpzc{M}_{5}&=&\{11111, 11121, 11211, 12111, 11221, 12211, 12121, 12221, 12321\}.
\end{eqnarray*}
Here are some examples of Motzkin words from $\mathpzc{AM}\setminus \mathpzc{M}$:
$$
2, 33, 232, 34433, 445654, 56565.
$$
Finally, the following words are not Motzkin words:
$$
12, 212, 23212, 23422, 2343212.
$$
In all these examples we used one-digit numbers to avoid confusion and keep the word notation, without commas or parentheses. 
If we want to use numbers that have more than one digit, we have to separate them with parentheses.
}
\end{Example}

\begin{figure}
\unitlength=1mm
\special{em.linewidth 0.5pt}
\linethickness{0.5pt}
\begin{picture}(140.00,110.00)(-45.00,30.00)

\put(15.00,115.00){\line(1,1){5.00}}
\put(20.00,120.00){\line(1,1){5.00}}
\put(25.00,125.00){\line(1,-1){5.00}}
\put(30.00,120.00){\line(1,-1){5.00}}
\put(15.00,115.00){\circle*{1.00}}
\put(20.00,120.00){\circle*{1.00}}
\put(25.00,125.00){\circle*{1.00}}
\put(30.00,120.00){\circle*{1.00}}
\put(35.00,115.00){\circle*{1.00}}
\put(20.50,112.00){\footnotesize $12321$}

\put(25.00,110.00){\line(0,-1){7.00}}

\put(15.00,95.00){\line(1,1){5.00}}
\put(20.00,100.00){\line(1,0){5.00}}
\put(25.00,100.00){\line(1,0){5.00}}
\put(30.00,100.00){\line(1,-1){5.00}}
\put(15.00,95.00){\circle*{1.00}}
\put(20.00,100.00){\circle*{1.00}}
\put(25.00,100.00){\circle*{1.00}}
\put(30.00,100.00){\circle*{1.00}}
\put(35.00,95.00){\circle*{1.00}}
\put(20.50,92.00){\footnotesize $12221$}

\put(25.00,90.00){\line(0,-1){7.00}}
\put(10.00,90.00){\line(-1,-1){7.00}}
\put(40.00,90.00){\line(1,-1){7.00}}

\put(15.00,75.00){\line(1,1){5.00}}
\put(20.00,80.00){\line(1,-1){5.00}}
\put(25.00,75.00){\line(1,1){5.00}}
\put(30.00,80.00){\line(1,-1){5.00}}
\put(15.00,75.00){\circle*{1.00}}
\put(20.00,80.00){\circle*{1.00}}
\put(25.00,75.00){\circle*{1.00}}
\put(30.00,80.00){\circle*{1.00}}
\put(35.00,75.00){\circle*{1.00}}
\put(20.50,71.00){\footnotesize $12121$}

\put(-15.00,75.00){\line(1,1){5.00}}
\put(-10.00,80.00){\line(1,0){5.00}}
\put(-05.00,80.00){\line(1,-1){5.00}}
\put(00.00,75.00){\line(1,0){5.00}}
\put(-15.00,75.00){\circle*{1.00}}
\put(-10.00,80.00){\circle*{1.00}}
\put(-05.00,80.00){\circle*{1.00}}
\put(00.00,75.00){\circle*{1.00}}
\put(05.00,75.00){\circle*{1.00}}
\put(-9.50,71.00){\footnotesize $12211$}

\put(-05.00,69.00){\line(0,-1){7.00}}
\put(03.00,69.00){\line(1,-1){7.00}}
\put(47.00,69.00){\line(-1,-1){7.00}}
\put(05.00,62.00){\line(1,1){7.00}}
\put(45.00,62.00){\line(-1,1){7.00}}

\put(45.00,75.00){\line(1,0){5.00}}
\put(50.00,75.00){\line(1,1){5.00}}
\put(55.00,80.00){\line(1,0){5.00}}
\put(60.00,80.00){\line(1,-1){5.00}}
\put(45.00,75.00){\circle*{1.00}}
\put(50.00,75.00){\circle*{1.00}}
\put(55.00,80.00){\circle*{1.00}}
\put(60.00,80.00){\circle*{1.00}}
\put(65.00,75.00){\circle*{1.00}}
\put(50.50,71.00){\footnotesize $11221$}

\put(55.00,69.00){\line(0,-1){7.00}}

\put(15.00,55.00){\line(1,0){5.00}}
\put(20.00,55.00){\line(1,1){5.00}}
\put(25.00,60.00){\line(1,-1){5.00}}
\put(30.00,55.00){\line(1,0){5.00}}
\put(15.00,55.00){\circle*{1.00}}
\put(20.00,55.00){\circle*{1.00}}
\put(25.00,60.00){\circle*{1.00}}
\put(30.00,55.00){\circle*{1.00}}
\put(35.00,55.00){\circle*{1.00}}
\put(20.50,51.00){\footnotesize $11211$}

\put(25.00,47.00){\line(0,-1){7.00}}
\put(03.00,45.00){\line(1,-1){7.00}}
\put(47.00,45.00){\line(-1,-1){7.00}}

\put(-15.00,55.00){\line(1,1){5.00}}
\put(-10.00,60.00){\line(1,-1){5.00}}
\put(-05.00,55.00){\line(1,0){5.00}}
\put(00.00,55.00){\line(1,0){5.00}}
\put(-15.00,55.00){\circle*{1.00}}
\put(-10.00,60.00){\circle*{1.00}}
\put(-05.00,55.00){\circle*{1.00}}
\put(00.00,55.00){\circle*{1.00}}
\put(05.00,55.00){\circle*{1.00}}
\put(-9.50,51.00){\footnotesize $12111$}

\put(45.00,55.00){\line(1,0){5.00}}
\put(50.00,55.00){\line(1,0){5.00}}
\put(55.00,55.00){\line(1,1){5.00}}
\put(60.00,60.00){\line(1,-1){5.00}}
\put(45.00,55.00){\circle*{1.00}}
\put(50.00,55.00){\circle*{1.00}}
\put(55.00,55.00){\circle*{1.00}}
\put(60.00,60.00){\circle*{1.00}}
\put(65.00,55.00){\circle*{1.00}}
\put(50.50,51.00){\footnotesize $11121$}

\put(15.00,35.00){\line(1,0){5.00}}
\put(20.00,35.00){\line(1,0){5.00}}
\put(25.00,35.00){\line(1,0){5.00}}
\put(30.00,35.00){\line(1,0){5.00}}
\put(15.00,35.00){\circle*{1.00}}
\put(20.00,35.00){\circle*{1.00}}
\put(25.00,35.00){\circle*{1.00}}
\put(30.00,35.00){\circle*{1.00}}
\put(35.00,35.00){\circle*{1.00}}
\put(20.50,31.00){\footnotesize $11111$}

\end{picture}
\caption{Lattice of Motzkin paths $\mathpzc{M}_{5}$ and associated reduced Motzkin words.
The greatest element of $\mathpzc{M}_{5}$ is associated with the word $12321$ and the 
least element is associated with the word $1^5$.}
\end{figure}
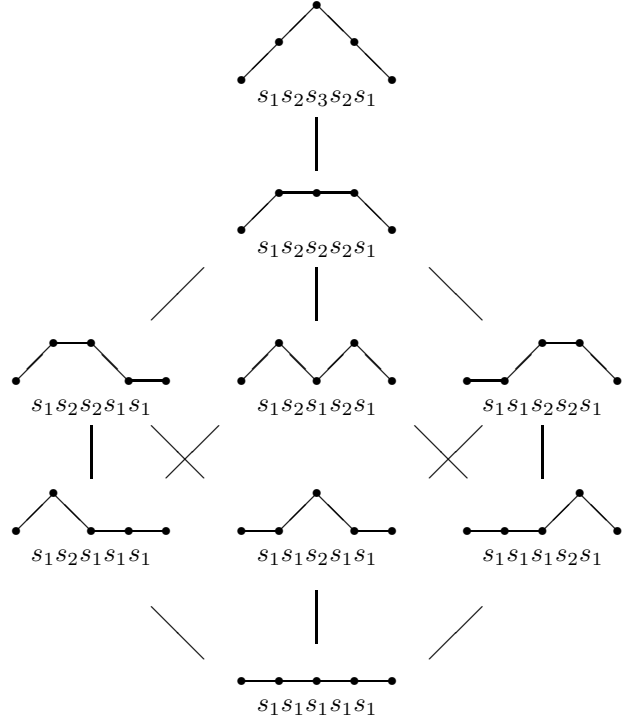

\begin{Remark}
{\rm 
Let us make a few remarks on Motzkin words.
\begin{enumerate}[(a)]
\item
Note that each set $\mathpzc{M}_{n}$ contains one word of the form $1^{n}$. These
words will be called {\it constant} since they correspond to constant Motzkin paths. 
Similarly, Motzkin words of the form $j^{n}$ will be called {\it constant}. 
The number of reduced Motzkin words of lenght $n$, given by the so-called Motzkin number $M_{n-1}$, 
(we need to use $n-1$ since our words are longer by one than the standard Motzkin words)
grows rapidly as $n$ increases. If we ignore the empty word and begin with words of length one, 
we get the following sequence: 1, 1, 2, 4, 9, 21, 51, 127, etc.
\item
With each reduced Motzkin word $w=j_1j_2\ldots j_n$ of positive lenght 
one can associate a {\it Motzkin path}, which is the graph of a nonnegative continuous 
function $f_w$ defined on the interval $[0,n-1]$ which starts from 
the origin $(0,0)$, ends at $(n-1,0)$, and consists of segments of three types: horizontal segments of lenght 1, 
NE segments from $(p,q)$ to $(p+1,q+1)$ and SE segments from $(p,q)$ to $(p+1,q-1)$. A bijection between reduced 
Motzkin words and Motzkin paths is given by the mapping $w\rightarrow f_w$, where $f_w(i)=j_i-1$. 
We understand that there is one Motzkin path corresponding to $n=1$, namely the one that reduces to one point.
In our paper, Motzkin words of height greater than one will be subwords of reduced Motzkin words and that is why
we can associate with them subpaths of Motzkin paths.
\item
It is well known that for each $n\geq 1$ the set of Motzkin paths of lenght $n$, and thus $\mathpzc{M}_{n}$, 
is a lattice, with partial order given by $f_1\leq f_2$ if and only if $f_1(x)\leq f_2(x)$ for any $x\in [0,n]$.
The maximal element of $\mathpzc{M}_{n}$ is $1\cdots (k-1)k(k-1) \cdots 1$, when $n=2k-1$ is odd and
$1\cdots (k-1)kk(k-1) \cdots 1$ when $n-2k$ is even. The minimal element of 
$\mathpzc{M}_{n}$ is the constant word of the form $1^{n}$.
\item
We do not use a standard terminology for Motzkin words. 
In our language, a reduced Motzkin word is directly related to the heights of the Motzkin paths at consecutive 
integers. Namely, if the height of a Motzkin path given by $f_w$ is $p$ at $i\in [k-1]$, then 
the corresponding letter from our alphabet $\mathbb{N}$ has color $p+1$.
Typically, reduced Motzkin words are built from letters which do not correspond to the heights of Motzkin 
paths, but to segments from which the Motzkin path is constructed. In that setting, if we assign letters $H,U,D$ to the horizontal, NE and SE segments, respectively, each reduced Motzkin word becomes shorter by one than in our 
notation. Using this notation, our reduced Motzkin words of lenght $k=4$ take the form
$$
HHH,\; HUD, \;UDH, \;UHD,
$$
respectively. For instance, our $1221$ corresponds to $UHD$ 
(1st segment goes up, 2nd segment is horizontal, 3rd segment goes down).
For our purposes, the standard language is inconvenient since to each word $w$ of lenght $n$ 
we assign a tuple of $n$ random variables and each letter $j$ encodes 
the information that this random variable is the $j$th replica of some indeterminate. 
\item
If $w=j_1j_2\cdots j_n$ is a Motzkin word (or, more generally, a word), then a {\it subword}
of $w$ is any word of the form $v=j_{i(1)}\ldots j_{i(p)}$ for some $i(1),\ldots, i(p)\in [n]$, where
$i(1)<\cdots <i(p)$. Thus, for instance, if 
$$
w=12^3323^24321,
$$
then the following words are subwords of $w$: $w_1=123$, $w_2=1243$, $w_3=323^24$, $w_4=2^332$,
but only the last one is a Motzkin word and only the last two are words obtained by taking consecutive letters in $w$.
Of course, a given subword of $w$, like $32$ in the considered example, may appear more than once in a given word. 
However, when speaking of a subword $v=j_{i(1)}\cdots j_{i(p)}$, we will usually 
understand that it encodes not only the information about the colors of letters, namely $j_{i(1)}, \ldots , j_{i(p)}$, 
but also the information about the corresponding indices $i(1), \ldots, i(p)$ (positions in $w$). 
\end{enumerate}
}
\end{Remark}

\section{Motzkin functionals}

We return to the study of mixed moments of the orthogonal replicas. It turns out that 
they reduce to the mixed moments associated with Motzkin words. All reduced 
Motzkin words correspond to moments of free random variables, whereas constant Motzkin words 
correspond to moments of boolean random variables.

\begin{Proposition}
Let $(\mathcal{A}_{i}, \varphi_{i})_{i\in I}$ be a family of noncommutative probability spaces and 
let $(\mathcal{A}_{{\rm rep}},\Phi)$ be the associated replica space. 
If $a_k\in {\mathcal A}_{i_k}$, where $k=1, \ldots , n$, then
$$
\Phi\left(a_1(j_1)\cdots a_{n}(j_n)\right)=0
$$ 
whenever $w:=j_1\cdots j_n\notin \mathpzc{M}_{n}(i_1, \ldots, i_n)$, where
$$
\mathpzc{M}_{n}(i_1, \ldots, i_n)=\{w\in \mathpzc{M}_{n}: j_k=j_{k+1} \;if\;i_k=i_{k+1}, \; k\in [n-1]\}.
$$
In particular, if $i_1\neq \ldots \neq i_n$, then the above equation holds whenever $w\notin \mathpzc{M}_{n}$.
\end{Proposition}
{\it Proof.}
If $i_1, \ldots , i_n$ are arbitrary, then the above moment vanishes in the 
following three cases:
\begin{enumerate}
\item $j_1\neq 1$ or $j_n\neq 1$, by the reduction property of Lemma 3.1,
\item $i_k=i_{k+1}$ and $j_{k}\neq j_{k+1}$ for some $k$, by orthogonality of replicas with the same labels
and different colors, by Proposition 3.2(1),
\item $i_k\neq i_{k+1}$ for some $k$ and $|j_{k}-j_{k+1}|\notin \{0,1\}$, by a repeated use of both the freeness and the reduction properties of Lemma 3.1 as well as the homomorphic property of replicas, namely
$$
a(j)a'(j)=(aa')(j)
$$
for any $a,a'\in \mathcal{A}_{i}$ and any $i,j$.
\end{enumerate}
Let us justify the last statement in more detail.
Namely, the freeness and reduction properties of Lemma 3.1 imply that nonzero moments of replicas must have 
consecutiveness of their colors for replicas with different labels, whereas the homomorphic property of Lemma 3.1 
together with the orthogonality of Proposition 3.2 imply that they must have stable colors for replicas with the 
same labels. As for the case when we use the homomorphic property, elements with the same labels 
meet and their product produces (apart from a kernel part) a constant multiple of the corresponding unit, which in fact 
leads to a deletion of this unit. Since originally it had a color bigger by one than its predecessor, 
this essentially means that the part with a constant multiple of the unit must have 
a color bigger by one than its successor or equal to it for the product to be nonzero. 
Note that if we are left at the end with a product of elements 
from the kernels and with different labels, we get zero contribution form 
such a moment. Therefore, 
the only situation when we get a nonzero contribution is when the sequence of colors $(j_1, \ldots, j_n)$ ends with
$j_n=1$. All these facts imply that the moment vanishes unless $w\in \mathpzc{M}_{m}(i_1, \ldots , i_n)$. 
If $i_1\neq \ldots \neq i_n$, then $\mathpzc{M}_{m}(i_1, \ldots , i_n)=\mathpzc{M}_{m}$ and thus the second 
property also holds. The proof is completed.
\hfill $\blacksquare$\\

Once we have the Motzkin structure underlying mixed moments of replicas, we can 
reformulate the results of Lemma 3.1, using a more elegant dual language, by defining
a family of path-dependent functionals on the free product of algebras (this can also be done for 
*-algebras). In that language, instead of 
looking simultaneously at variables associated with different notions of independence, we look 
simultaneously at different ways of convolving these variables by means of product functionals.

Therefore, let us define a general family of multilinear functionals on the free product of unital algebras ${\mathcal A}_{i}$, where 
$i\in I$ (the unit in ${\mathcal A}_{i}$ will be denoted by $1_{i}$), in which we do not identify units, 
$
{\mathcal A}:=*_{i\in I}\mathcal{A}_{i},
$
with the unit $1_{{\mathcal A}}$. Instead of a sequence of functionals, 
we will introduce a family of {\it path-dependent functionals}. We allow their arguments to 
be arbitrary polynomials, say $a_k=f(b_1, \ldots , b_m)$, where $f$ is a polynomial in noncommutative 
variables $b_1\in {\mathcal A}_{i_1}, \ldots , b_m\in \mathcal{A}_{i_m}$.
The corresponding replicas will be defined by $a_k(j):=f(b_1(j), \ldots , b_m(j))$ for any $j\in \mathbb{N}$. 
Nevertheless, as in the case of sequences of moment functionals, it will be sufficient 
to assume that each $a_k$ is a monomial from some ${\mathcal A}_{i_k}$ or the unit 
$1_{{\mathcal A}}$.

\begin{Definition}
{\rm 
With the notations of Proposition 5.1, let ${\mathcal A}:=*_{i\in I}\mathcal{A}_{i}$ be the 
free product with non-identified units and with the unit $1_{{\mathcal A}}$. Define a family of 
path-dependent multilinear functionals 
$$
\{\varphi(w):w\in \mathpzc{M}^{*}\}
$$ 
on ${\mathcal A}^{|w|}$ as multilinear extensions of
\begin{eqnarray*}
\varphi(\varnothing)(1_{{\mathcal A}})&=&1\\
\varphi(w)\left(a_1,\ldots ,a_n\right):&=&\Phi\left(a_{1}(j_1)\cdots a_{n}(j_n)\right)
\end{eqnarray*}
for any $a_k\in \mathcal{A}$, $k=1, \ldots , n$, where $w=j_1\cdots j_n$.
The functionals $\varphi(w)$ will be called {\it multilinear Motzkin functionals}. 
}
\end{Definition}

\begin{Lemma}
Let $w=j_1\cdots j_n\in \mathpzc{M}_{n}$ and let $a_k\in \mathcal{A}_{i_k}$.
Then we have
\begin{enumerate}
\item
the freeness property
$$
\varphi(w)\left(a_1,\ldots ,a_n\right)=0
$$
whenever $a_k\in {\rm Ker}(\varphi_{i_k})$ for $k=1, \ldots ,n$, where $i_1\neq \cdots \neq i_n$,
\item the reduction property
$$
\varphi(w)\left(a_1,\ldots, 1_{i(r)},\ldots ,a_n\right)=\varphi(\check{w})\left(a_1,\ldots ,\check{1}_{i(r)}, \ldots ,a_n\right)
$$
whenever $a_k\in {\rm Ker}(\varphi_{i_k})$ for 
$k=1, \ldots ,r-1$ and $(j_1, \ldots , j_r)=(1,\ldots , r)$, where $i_1\neq \cdots \neq i_r$ and
$\check{w}=j_1\cdots j_{r-1}j_{r+1}\cdots j_n$ and for other colors this moment vanishes,
\item the homomorphic property
$$
\varphi(w)\left(\ldots, a_{k}, a_{k+1}, \ldots\right)=\varphi(\check{w})\left(\ldots, a_{k}a_{k+1}, \ldots\right)
$$
whenever $(i_k,j_k)=(i_{k+1},j_{k+1})$, where $\check{w}=j_1\cdots j_{k}j_{k+2}\cdots j_n$ 
and this moment vanishes if $i_k=i_{k+1}$ and $j_k\neq j_{k+1}$.
\end{enumerate}
\end{Lemma}
{\it Proof.}
Properties (1) and (2) are straightforward reformulations of properties of Lemma 3.1. Property (3) 
follows from the definition of replicas and their orthogonality.
\hfill $\blacksquare$\\

Lemma 5.1 gives the main ingredients that lead to the decomposition of the free product
of functionals. An important feature of our theory is that it gives a `refinement' of the free product of functionals 
(states, in the case of *-noncommutative probability spaces) in which the functionals 
(states) $\varphi(w)$ play the role of a generating set of a set of moment functionals.
These functionals are multilinear, which is natural in the context of Motzkin words, but 
of course one can also define related linear moment functionals $\psi(w)$ in a natural way, which will be done below.
The whole family $\{\psi(w):w\in \mathpzc{M}^*\}$ will be needed to describe the free product of functionals.
Moreover, the boolean product of states is included in this framework: 
it suffices to take $\{\psi(w):w\in \mathpzc{B}^{*}\}$, where 
$\mathpzc{B}^{*}=\{1^{n}, n\in \mathbb{N}\}$, which corresponds to constant Motzkin paths. This will 
be the subject of Theorem 5.1 below.

\begin{Remark}
{\rm 
In order to state our results on the decompositions of the free and boolean products of functionals 
in terms of Motzkin functionals, we would like to use linear moment functionals. Let us recall some 
basic facts.
\begin{enumerate}[(a)]
\item
The free product of algebras ${\mathcal A}_{i}$ 
with identification of units will be denoted by 
$$
\star_{i\in I}{\mathcal A}_{i}
$$
and it can be viewed as the quotient algebra ${\mathcal A}/U$, where 
${\mathcal A}$ is the free product without identification of units and 
$U=\langle \{1_{{\mathcal A}}-1_i: i\in I\}\rangle$, the two-sided ideal generated by all $1_{{\mathcal A}}-1_i$. 
The free product of functionals $\varphi_i$ on $\star_{i\in I}{\mathcal A}_{i}$
in the sense of Voiculescu will be denoted by $\star_{i\in I}\varphi_i$ (we use the notation of Nica \cite{[19]}).
\item
Decompose ${\mathcal A}$ into a vector space direct sum 
$$
{\mathcal A}:=\bigoplus_{n=0}^{\infty} F_{n},
$$
where $F_n$ is spanned by all monomials $a_{1}\cdots a_{n}$ of degree $n$, where $a_k\in {\mathcal A}_{i_k}$ 
for $k=1, \ldots , n$ and $i_1\neq \cdots \neq i_n$
and $F_{0}={\mathbb C}1_{{\mathcal A}}$. 
An arbitrary noncommutative polynomial $f\in {\mathcal A}$ can be uniquely decomposed as 
a sum of polynomials $f_n$, namely
$$
f=\sum_{n=0}^{\infty}f_n,
$$
where $f_n\in F_n$ for all $n$ and only a finite number of terms is nonzero (uniqueness refers to the direct sum decomposition
since within each $F_n$ we may have different realizations of the same element since algebras ${\mathcal A}_{i}$ are arbitrary).
\end{enumerate}
}
\end{Remark}

\begin{Definition}
{\rm Let $\Phi$ be the tensor product state given by Definition 3.1. 
Define a family of path-dependent linear functionals
$$
\mathpzc{G}(\{\varphi_i:i\in I\}):=\{\psi(w):w\in \mathpzc{M}^{*}\},
$$
where $\psi(w):{\mathcal A}\rightarrow \mathbb{C}$ are given by linear extensions of
\begin {eqnarray*}
\psi(\varnothing)&:=&Id_{\,{\mathbb C}1}\\
\psi(w)(a_{i_1} \cdots  a_{i_n})&:=&\Phi(a_{i_1}(j_1) \cdots a_{i_n}(j_n))
\end{eqnarray*}
for $w=j_1\cdots j_n$, where $a_k\in \mathcal{A}_{i_{k}}$ for any $k$ and 
$i_1\neq \cdots \neq i_n$, $n\in {\mathbb N}$, and are set to be zero on $F_m$ for $m\neq n$.
The functionals $\psi(w)$ will be called {\it linear Motzkin functionals}.
}
\end{Definition}

\begin{Remark}
{\rm 
We would like to make some comments on $\varphi(w)$ and $\psi(w)$.
\begin{enumerate}[(a)]
\item
Multilinear and linear Motzkin functionals are closely related to each other and they determine each other uniquely.
In particular, $$\varphi(w)(a_1, \ldots, a_n)=\psi(w)(a_1\cdots a_n)$$ whenever $a_k\in {\mathcal A}_{i_k}$ 
and $i_1\neq \ldots \neq i_n$ and in that situation it is irrelevant 
which functionals we use. However, when some of the neighboring indices are the same, then 
quite often the multilinear functionals seem to be more appropriate since they respect the order of 
each monomial. For instance, the $n$th moment of a replica $a(1)$, where $a\in \mathcal{A}_{i}$, 
can be described by
$$
\varphi(1^n)(a,\ldots , a)=\Phi(a(1)^{n})\;\;\;{\rm and}\;\;\psi(1)(a^n)=\Phi(a^n(1))
$$
and in some context, like that of convolutions and cumulants, the first description is more convenient 
since $w=1^n\in \mathpzc{M}_{n}$ is used to encode the information on moments of order $n$.
On this occasion, note that the freeness and reduction properties have their counterparts for linear Motzkin functionals.
\item
We are ready to use a `duality trick' that allows us to convolve variables living in the free product of algebras
in many different ways simultaneously. For any polynomial $f\in F_n$, the following {\it duality relation} holds:
$$
\psi(w)(f)=\Phi(f(w)),
$$
where $f(w)$ is a linear combination of products of $n$ elements from different algebras (with different neighbors). 
In this fashion, 
elements of the lattice ${\mathpzc M}_n$, which are originally arguments of random variables, become
arguments of functionals on the algebra generated by these variables. 
In this connection, observe that we can take these arguments to be elements of the semigroup
algebra ${\mathbb C}(\mathpzc{M}^{*})$ by linearity:
$$
\psi(\alpha_1 w_1+\alpha_2 w_2):=\alpha_1 \psi(w_1)+\alpha_2 \psi(w_2)
$$
for any $w_1,w_2\in {\mathpzc M}^{*}$ and any $\alpha_1, \alpha_2 \in \mathbb{C}$.
Thus, to the linearity at the level of variables we are adding linearity and at the level of functionals, 
controlled by ${\mathbb C}(\mathpzc{M}^{*})$.
\end{enumerate}
}
\end{Remark}

\begin{Definition}
{\rm 
More generally, let $\gamma =\sum_{w\in \mathpzc{M}^{*}}\gamma_{w}w\in {\mathbb C}\langle \langle \mathpzc{M}^{*}\rangle \rangle$, 
understood as the formal power series on $\mathpzc{M}^{*}$.
Define the linear functionals on ${\mathcal A}$ by
$$
\psi(\gamma):=\sum_{w\in \mathpzc{M}^{*}}\gamma_{w}\psi(w),
$$
where, by Definition 5.1, only a finite number of terms do not vanish when acting on ${\mathcal A}$.
The linear space
$$
\mathpzc{M}\left(\{\varphi_i:i\in I\}\right):=\{\psi(\gamma): \,\gamma\in {\mathbb C}\langle\langle \mathpzc{M}^{*}\rangle \rangle\}
$$
will be called the {\it space of Motzkin product functionals} of $\varphi_i$. If $\gamma_{\emptyset}=1$, the 
associated functional is normalized. The family $\mathpzc{G}\left(\{\varphi_i:i\in I\}\right)$ 
plays the role of its generating set. 
}
\end{Definition}

Using this generating set, we can decompose the free product of normalized linear functionals 
of Avitzour \cite{[2]} and Voiculescu \cite{[25]}. On this occasion, we also decompose the boolean product 
of normalized linear functionals. Equivalently, we can represent these
products as Motzkin product functionals (modulo the identification of units in the case of the free product).

\begin{Theorem}
If $\{(\mathcal{A}_{i}, \varphi_i):i\in  I\}$ is a family of noncommutative probability spaces, then
\begin{enumerate}
\item
for the free product of functionals $\varphi_i$, it holds that
\begin{eqnarray*}
(\star_{i\in I}\varphi_{i})\circ \tau &=&\sum_{w\in \mathpzc{M}^{*}}\psi(w),
\end{eqnarray*}
where $\tau: {\mathcal A}\rightarrow \star_{i\in I}{\mathcal A}_{i}$ is the canonical unital
homomorphism that maps each internal unit $1_i$ onto $1$,
\item
for the boolean product $\convolution_{i\in I}\varphi_i$ of functionals, it holds that
\begin{eqnarray*}
\convolution_{i\in I}\varphi_i&=&\sum_{w\in \mathpzc{B}^{*}}\psi(w),
\end{eqnarray*}
where $\mathpzc{B}^{*}=\{1^{n}: n\in \mathbb{N}^{*}\}$ and $1^{0}=\varnothing$.
\end{enumerate}
\end{Theorem}
{\it Proof.}
It is easy to see that the $n$th moment of each $a\in {\mathcal A}_{i}$ with respect to $\varphi_i$ 
agrees with the moment of $a^n(1)$ with respect to $\Phi$ and thus it agrees with 
the $n$th moment of $a$ w.r.t. $\psi(1)$ for any $n$. Therefore, 
the distribution of $a$ with respect $\varphi_i$ agrees with its distribution 
under 
$$
\psi\big(\sum_{w\in \mathpzc{B}^*}w\big)\equiv\sum_{w\in \mathpzc{B}^*}\psi(w).
$$
Thus, it also agrees with its distribution
under 
$$
\psi\big(\sum_{w\in \mathpzc{M}^*}w\big)\equiv\sum_{w\in \mathpzc{M}^*}\psi(w)
$$ 
since $\psi(w)(a^n)=0$ if $w\neq 1$.
Now, let $i_1\neq \cdots \neq i_n$. 
Then
\begin{eqnarray*}
\psi(1^n)(a_1\cdots a_n)&=&
\psi(1^n)(a_1^{\circ}a_2 \cdots a_n)+
\varphi_{i_1}(a_1)\psi(1^{n})(1_{i_1}a_2\cdots a_n)\\
&=& \varphi_{i_1}(a_1)\psi(1^{n-1})(a_2\cdots a_n)
\end{eqnarray*}
by the reduction property, where $a_1^{\circ}=a_1-\varphi_{i_1}(a_1)1_{i_1}$.
We can continue in this fashion to obtain the complete factorization which 
holds in the boolean case. Therefore,  
$$
\convolution_{i\in I}\varphi_i=\sum_{w\in \mathpzc{B}^{*}}\psi(w).
$$
In turn, by the freeness property of Lemma 5.1,
$$
\sum_{w\in \mathpzc{M}^{*}}\psi(w)(a_1^{\circ}\cdots a_n^{\circ})=0,
$$
where $a_k^{\circ}=a_{k}-\varphi_{i_k}(a_{k})1_{i_k}$ for all $k=1, \ldots , n$,
and thus $\psi\left(\sum_{w\in \mathpzc{M}^{*}}w\right)$ agrees with 
$(\star_{i\in I}\varphi_{i})\circ \tau$ on such products, where $\tau(1_{i})=1$, the unit
in $\star_{i\in I}\mathcal{A}_{i}$, for all $i$.
We still need to justify that
$$
\sum_{w\in \mathpzc{M}_{n}}\psi(w)(a_1^{\circ}\cdots a_{r-1}^{\circ}1_{i_r}a_{r+1}\cdots a_n)=
\sum_{w\in \mathpzc{M}_{n-1}}\psi(w)(a_1^{\circ}\cdots a_{r-1}^{\circ}a_{r+1}\cdots a_n)
$$
for any $r\geq 1$. By the implementation of the reduction property to 
linear functionals $\psi(w)$, we have 
\begin{eqnarray*}
LHS&=&\sum_{w=1\cdots rw'\in \mathpzc{M}_{n}}\psi(w)(a_1^{\circ}\cdots a_{r-1}^{\circ}1_{i_r}a_{r+1}\cdots a_n)\\
&=& \sum_{w=1\cdots (r-1)w'\in \mathpzc{M}_{n-1}}\psi(w)(a_1^{\circ}\cdots a_{r-1}^{\circ}a_{r+1}\cdots a_n).
\end{eqnarray*}
The last equality holds since if the $(r+1)$th letter is equal to $r+1$, then after deleting $r$ we get the word
$w=1\cdots (r-1)w' \notin \mathpzc{M}_{n-1}$, where $w'$ begins with $r+1$, and the contribution from such $\psi(w)$
is zero (if $a_{r+1}\in {\rm Ker}(\varphi_{i_{r+1}})$, this follows from the freeness property, and if $a_{r+1}=1_{i_{r+1}}$, 
we use the reduction property).
Therefore, under $\psi\left(\sum_{w\in \mathpzc{M}^{*}}w\right)$ all units $1_{i}$ are identified with $1_{{\mathcal A}}$
and thus the canonical homomorphism $\tau$ given by $\tau(a)=a+I$, where 
$I$ is the two-sided ideal generated by the set $\{1_{i}-1_{{\mathcal A}}: i\in I\}$
and $1=1_{{\mathcal A}}+I$ is the unit in $\star_{i\in I}\varphi_{i}$, 
satisfies the desired relation. In other words,  $\psi(\sum_{w\in \mathpzc{M}^{*}}w)$ factors 
through $\star_{i\in I}\varphi_i$ and the proof is completed.
\hfill $\blacksquare$\\

Using the language of mixed moments, Theorem 5.1(1) gives a nice formula for the mixed moments of 
free random variables under the free product of normalized linear functionals 
$\star_{i\in I}\varphi_i:\star_{i\in I}\mathcal{A}_{i}\rightarrow {\mathbb C}$.

\begin{Corollary}
The mixed moments of free random variables can be expressed in the form
$$
\star_{i\in I}\varphi_i (a_1\cdots a_n)=\sum_{w=j_1\cdots j_n\in \mathpzc{M}_n}\Phi(a_1(j_1)\cdots a_n(j_n))
$$
for any $a_1\in \mathcal{A}_{i_1}, \ldots, a_n\in \mathcal{A}_{i_n}$, where $i_1\neq \cdots \neq i_n$ and 
$a_k(j_k)$ is the orthogonal replica of $a_k$ of color $j_k$ for any $k=1, \ldots, n$. 
\end{Corollary}

Let us summarize the above results: 
\begin{enumerate}
\item
Proposition 5.1 exhibits the Motzkin structure on the level of the mixed moments of replicas,
\item
Lemma 5.1 establishes conditions on the Motzkin functionals which give a `refinement' of freeness,
\item
Theorem 5.2 gives the decomposition of the free and boolean products of functionals in terms of Motzkin functionals.
\end{enumerate}

It is not hard to see that the Motzkin path decomposition of the 
conditionally free product of functionals is easy to obtain: it suffices to modify $\Psi$ and construct it from pairs of
functionals on each ${\mathcal A}_{i}$ as it was done in \cite{[11]} (in fact, even a countable
number of functionals for each algebra causes no problems).
We have also obtained similar results for nonsymmetric products like the orthogonal, s-free and monotone products, 
in the case of two noncommutative probability spaces.
For that purpose, we need to repeat what we have done above, except that we need to take
smaller families of replicas, as shown in Theorem 3.2. Details will be given in a separate paper.
 
The exact values of the moments of replicas under Motzkin functionals, expressed in terms of boolean cumulants
and some noncrossing partitions suitably adapted to $w\in \mathpzc{M}$ will be given in Section 8.

\section{Moments}

Let us present some useful formulas for moment functionals $\psi(w)$. 
They are recursions, as that in Lemma 5.1, but they are more convenient since 
they involve arbitrary variables $a_1, \ldots, a_n$ rather than those in the kernels
of the corresponding functionals. 

The first formula reminds the one of Avitzour for free random variables \cite{[2]}, where 
singletons are also pulled out, but the remaining elements are from the kernels of the corresponding states. 
We denote by $\mathpzc{A}$ the set of artihmetic progressions of the form $\mathpzc{A}=\{1, 12, 123, \ldots\}$. 
In particular, if we write that some sequence belongs to $\mathpzc{A}$, we understand that it is non-empty.

\begin{Lemma}
Let $a_k\in \mathcal{A}_{i_k}$, where $i_1\neq \cdots \neq i_n$ and $k=1, \ldots, n$.
\begin{enumerate}
\item
For any $w\in \mathpzc{M}_n$, it holds that
$$
\psi(w)\left(a_1\cdots a_n\right)=
\sum_{\stackrel{[n]=A\cup B}{\scriptscriptstyle j_1\cdots  j_{{\rm max}A}}\in \mathpzc{A}}
(-1)^{|A|-1}\prod_{k\in A}\varphi_{i_k}(a_{k})
\cdot \psi(w_{B})(a_B),
$$
where products $w_{B}:=\prod_{k\in B}j_k$ and
$a_B:=\prod_{k\in B}a_{k}$ are taken in the order of increasing $k$.
\item
In particular, if $w=1w'$, where $w'\in \mathpzc{M}_{n-1}$, then 
$$
\psi(w)\left(a_1\cdots a_n\right)=\varphi_{i_1}(a_1)\psi(w')(a_2\cdots a_n)
$$
for $n>1$. 
\item
Thus, if $w=1^n$, then 
$$
\psi(w)\left(a_1\cdots a_n\right)=\varphi_{i_1}(a_1)\cdots \varphi_{i_n}(a_n),
$$
which gives the moments of the boolean product $\convolution_{i\in I} \varphi_i$.
\end{enumerate}
\end{Lemma}
{\it Proof.}
It suffices to prove (1) since (2) and (3) are easy consequences of (1). We shall 
prove (1) by induction, using typical decompositions $a_k=a_k^{\circ}+\varphi_{i_k}(a_1)1_{i_k}$
for $k=1, \ldots, n$. 
We claim that for any $k=1, \ldots ,n$ the following formula holds:
\begin{eqnarray*}
\psi(w)(a_1 \cdots a_n)&=&\psi(w)(a_1^{\circ}\cdots a_k^{\circ}a_{k+1}\cdots, a_n)\\
&+&
\sum_{\stackrel{[n]=A\cup B}{\scriptscriptstyle {\rm max}A \leq k}}
(-1)^{|A|-1}\prod_{k\in A}\varphi_{i_k}(a_{k})
\cdot \psi(w_{B})(a_{B})\delta_{j_1,1}\cdots \delta_{j_k,k}.
\end{eqnarray*}
We will prove it by induction.
The first step gives
$$
\psi(w)(a_1\cdots a_n)=\psi(w)(a_1^{\circ} a_2 \cdots a_n)+
\varphi_{i_1}(a_1)\psi(\check{j_1}\cdots j_n)(\check{a}_1 \cdots a_n)\delta_{j_1,1}
$$
by the reduction property.
Suppose that the formula holds for some $k$. We will prove that it holds for $k+1$.
Using the decomposition $a_{k+1}=a_{k+1}^{\circ}+\varphi_{i_{k+1}}(a_{k+1})1_{i_{k+1}}$ in the first term on the RHS
of the above formula, we obtain
\begin{eqnarray*}
\psi(w)(a_1^{\circ}\cdots a_k^{\circ}a_{k+1} \cdots a_n)
&=&
\psi(w)(a_1^{\circ}\cdots a_{k+1}^{\circ} a_{k+2}\cdots a_n)
\\
&+&
\varphi_{i_{k+1}}(a_{k+1})\psi(w)(a_1^{\circ} \cdots a_k^{\circ} 1_{i_{k+1}} \cdots a_n)\delta_{j_{k+1},k+1}
\end{eqnarray*}
The first term is of the desired form. Writing $a_j^{\circ}=a_{j}-\varphi_{i_j}(a_j)1_{i_j}$ 
for all $j=k, \ldots, 1$ in the second term (in that order), we obtain 
$$
\varphi_{i_{k+1}}(a_{k+1})\psi(w)(a_1^{\circ} \cdots a_{k}^{\circ} 1_{i_{k+1}}\cdots a_n)1_{i_{k+1},k+1}
$$
$$
=
\sum_{\stackrel{[n]=A\cup B}{\scriptscriptstyle {\rm max}A = k+1}}
(-1)^{|A|-1}\prod_{k\in A}\varphi_{i_k}(a_{k})
\cdot \psi(w_{B})(a_{B})\delta_{j_1,1}\cdots \delta_{j_{k+1},k+1}
$$
since each subset $A'$ of $[k]$ associated with subtracted terms of the form
$\varphi_{i_j}(a_j)1_{i_j}$ produces $(-1)^{|A'|}\prod_{j\in A'}\varphi_{i_j}(a_j)$ and 
there is an additional singleton $\varphi_{i_{k+1}}(a_k)$, this gives 
$(-1)^{|A|-1}\prod_{k\in A}\varphi_{i_k}(a_{k})$ times $\psi(w_B)(a_B)$ by the reduction property used repeatedly.
This completes the proof of (1).

\hfill $\blacksquare$

\begin{Example}
{\rm
Let us compute the family of moments
$$
\{\psi(w)(a_1b_2a_3b_4):w\in \mathpzc{M}_{4}\},
$$
where $a_1,a_3\in \mathcal{A}_{1}$ and $b_2,b_4\in \mathcal{A}_{2}$.
Using Lemma 6.1(3), we obtain
\begin{eqnarray*}
\psi(1^4)(a_1b_2a_3b_4)&=&\varphi_1(a_1)\varphi_2(b_2)\varphi_1(a_3)\varphi_2(b_2)
\end{eqnarray*}
Next, using Lemma 6.1, first (2) and then (1), we obtain
\begin{eqnarray*}
\psi(1^221)(a_1b_2a_3b_4)&=&\varphi_1(a_1)\psi(121)(b_2a_3b_2)\\
&=&\varphi_1(a_1)\left(\varphi_1(a_3)\varphi_2(b_2b_4)-\varphi_1(a_3)\varphi_2(b_2)\varphi_2(b_4)\right),
\end{eqnarray*}
where the first term corresponds to the decomposition $[2,4]=\{3\}\cup \{2,4\}$ and the second one -- to 
$[2,4]=\{2,3\}\cup \{4\}$, with the convenient understanding that we decompose the interval $[2,4]=A\cup B$ rather than the set 
$[3]$, which is equivalent and allows us to keep track of the original indices.
Similarly, 
\begin{eqnarray*}
\psi(121^2)(a_1b_2a_3b_4)&=&\varphi_2(b_2)\psi(1^3)(a_1a_3b_4) - \varphi_1(a_1)\varphi_2(b_2)\psi(1^2)(a_3,b_4)\\
&=&\varphi_2(b_2)\left(\varphi_1(a_1a_3)\varphi_2(b_4)-\varphi_1(a_1)\varphi_1(a_3)\varphi_2(b_4)\right),
\end{eqnarray*}
where the first term corresponds to the decomposition $[1,4]=\{2\}\cup \{1,3,4\}$ and the second one -- to 
$[1,4]=\{1,2\}\cup \{3,4\}$. Finally, 
$$
\psi(12^21)(a_1b_2a_3b_4)=0
$$
since for any $A$, the word $w_{B}$ corresponding to the complement of $A$ is not a Moztkin word.
It is not hard to check that 
$$
\sum_{w\in \mathpzc{M}_{4}}\psi(w)(a_1b_2a_3b_4)=(\varphi_1\star \varphi_2)(a_1b_2a_3b_4),
$$
which agrees with Theorem 5.1. Computations are similar for higher order moments.
Let us observe that we always have to include only such decompositions $[n]=A\cup B$ which 
give consecutiveness of colors $j_1, \ldots ,j_{{\rm max} A}$. 
This limits the number of possibilities. 
}
\end{Example}

\begin{Corollary}
Under the assumptions of Lemma 6.2, the moments of free random variables can be expressed in the form
$$
\left(\left(\star_{i\in I}\varphi_i\right)\circ \tau \right)(a_1\cdots a_n)=
\sum_{w=j_1\cdots j_n\in \mathpzc{M}_{n}}
\sum_{\stackrel{[n]=A\cup B}{\scriptscriptstyle j_1\cdots  j_{{\rm max}A}}\in \mathpzc{A}}
(-1)^{|A|-1}\prod_{k\in A}\varphi_{i_k}(a_{k})\cdot \psi(w_B)\big(\prod_{k\in B}a_{k}\big),
$$
where the index $k$ is taken in the increasing order.
\end{Corollary}
{\it Proof.}
This is an immediate consequence of Theorem 5.1 and Lemma 6.1.
\hfill $\blacksquare$\\

\begin{figure}
\unitlength=1mm
\special{em.linewidth 0.5pt}
\linethickness{0.5pt}
\begin{picture}(140.00,50.00)(-27.00,-20.00)

\put(15.00,15.00){\line(1,1){5.00}}
\put(20.00,20.00){\line(1,1){5.00}}
\put(25.00,25.00){\line(1,-1){5.00}}
\put(30.00,20.00){\line(1,0){5.00}}
\put(35.00,20.00){\line(1,-1){5.00}}
\put(40.00,15.00){\line(1,1){5.00}}
\put(45.00,20.00){\line(1,-1){5.00}}
\put(50.00,15.00){\line(1,1){5.00}}
\put(55.00,20.00){\line(1,1){5.00}}
\put(60.00,25.00){\line(1,-1){5.00}}
\put(65.00,20.00){\line(1,-1){5.00}}

\put(15.00,15.00){\circle*{1.00}}
\put(20.00,20.00){\circle*{1.00}}
\put(25.00,25.00){\circle*{1.00}}
\put(30.00,20.00){\circle*{1.00}}
\put(35.00,20.00){\circle*{1.00}}
\put(40.00,15.00){\circle*{1.00}}
\put(45.00,20.00){\circle*{1.00}}
\put(50.00,15.00){\circle*{1.00}}
\put(55.00,20.00){\circle*{1.00}}
\put(60.00,25.00){\circle*{1.00}}
\put(65.00,20.00){\circle*{1.00}}
\put(70.00,15.00){\circle*{1.00}}

\put(22.00,8.00){\footnotesize $w=   1232^21212321$}

\put(07.00,-10.00){\line(1,1){5.00}}
\put(17.00,-05.00){\line(1,1){5.00}}
\put(22.00,0.00){\line(1,-1){5.00}}
\put(27.00,-05.00){\line(1,0){5.00}}
\put(37.00,-05.00){\line(1,-1){5.00}}
\put(47.00,-10.00){\line(1,1){5.00}}
\put(52.00,-05.00){\line(1,-1){5.00}}
\put(57.00,-10.00){\line(1,1){5.00}}
\put(62.00,-05.00){\line(1,1){5.00}}
\put(67.00,0.00){\line(1,-1){5.00}}
\put(72.00,-05.00){\line(1,-1){5.00}}

\put(07.00,-10.00){\circle*{1.00}}
\put(12.00,-05.00){\circle*{1.00}}
\put(17.00,-05.00){\circle*{1.00}}
\put(22.00,0.00){\circle*{1.00}}
\put(27.00,-05.00){\circle*{1.00}}
\put(32.00,-05.00){\circle*{1.00}}
\put(37.00,-05.00){\circle*{1.00}}
\put(42.00,-10.00){\circle*{1.00}}
\put(47.00,-10.00){\circle*{1.00}}
\put(52.00,-05.00){\circle*{1.00}}
\put(57.00,-10.00){\circle*{1.00}}
\put(62.00,-05.00){\circle*{1.00}}
\put(67.00,0.00){\circle*{1.00}}
\put(72.00,-05.00){\circle*{1.00}}
\put(77.00,-10.00){\circle*{1.00}}

\put(16.00,-15.00){\footnotesize $w_1=232^2$}
\put(47.00,-15.00){\footnotesize $w_2=1212321$}

\end{picture}
\caption{Decomposition of a Motzkin path $w=1w_1w_2$.}
\end{figure}
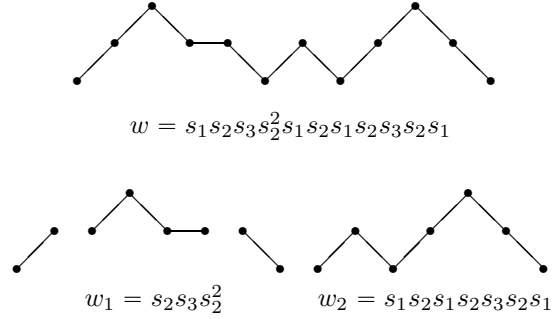

\begin{Remark}
{\rm 
Before we give another recursion for moments under Motzkin functionals, let us recall
the well-known recursion for Motzkin numbers:
$$
M_{n}=M_{n-1}+\sum_{k=1}^{n-2}M_kM_{n-k-1},
$$
where $M_1=1$ and $M_k=|\mathpzc{M}_{k}|$. Let us recall that our notation for Motzkin paths and thus for 
Motzkin numbers is slightly different than the standard one, thus $M_2=1, M_3=2, M_4=4, M_5=9$, etc.
A straightforward proof of this recursion is based on a decomposition of the paths from $\mathpzc{M}_{n}$ 
described in the proof of Lemma 6.2 given below. A simple illustration of such a decomposition is given in 
Fig. 2. The first return of the path to the $x$-axis divides the path into two subpaths, represented by 
words $1w_11\in \mathpzc{M}_{k+2}$ and $w_2\in \mathpzc{M}_{n-k-1}$. Then we remove the first and last segments
from the first subpath to get $w_1\in \mathpzc{AM}_{k}$ of height $2$. This path contributes to $M_k$ in the 
above recurrence since the cardinality of the set of Motzkin subpaths corresponding to Motzkin words 
of height $2$ and lenght $k$ is clearly equal to $M_k$.
}
\end{Remark}

\begin{Lemma}
Let $a_j\in \mathcal{A}_{i_j}$, where $i_1\neq \cdots \neq i_n$.
If $w\in \mathpzc{M}_{n}$ is of the form 
$w=1w_1w_2$, where $w_1\in \mathpzc{AM}_{k}$ is of height $h(w_1)=2$, and $w_2\in \mathpzc{M}_{n-k-1}$, then
\begin{eqnarray*}
\psi(w)\left(a_1\cdots a_n\right)&=&\psi(\widetilde{w}_1)(a_2 \cdots a_{k+1})\psi(1w_2)(a_1a_{k+2}\cdots a_n)\\
&-&\varphi_{i_1}(a_1)\psi(\widetilde{w}_1)(a_2\cdots a_{k+1})\psi(w_2)(a_{k+2} \cdots a_{n}),
\end{eqnarray*}
where $\widetilde{w}_{1}$ is the reduced word obtained from $w_1$ by mapping each letter $j$ onto $j-1$.
\end{Lemma}
{\it Proof.}
The path that corresponds to $w\in \mathpzc{M}_{n}$ begins with a segment $(1,1)$ and it can be
decomposed into two parts: the subpath starting at $(0,0)$ and ending at the first return of the 
path $w$ to the $x$-axis at $(k,0)$ for some $k>1$, and the remaining subpath starting 
from $(k,0)$ and ending at $(n,0)$, with possible returns to the $x$-axis.
The words corresponding to these subpaths are of the form
$$
1w_11\in \mathpzc{M}_{k+2}\;\;\;{\rm and}\;\;\;w_{2}\in \mathpzc{M}_{n-k-1},
$$
where $w_1$ is a Motzkin word of height $h(w_1)=2$. Note that the first 
subpath corresponds to the sequence of first $k+2$ arguments of $\varphi(w)$, 
namely $(a_{1}, \ldots , a_{k+2})$. Now, replica $a_{k+2}(1)$ multiplies by a cyclic projection from the right all 
elements at tensor sites of all colors and labels $\neq i_{k+2}$. 
In turn, the elements with label $i_{k+2}$ and color $>1$ are multiplied from the right by a 
cyclic projection in $a_{k+1}(2)$ for $k>1$ since $i_{k+1}\neq i_{k+2}$ 
(if $k\leq 1$, there are no such elements).
This fact allows us to pull out the moment of a Motzkin functional corresponding to the word 
$w_2\in \mathpzc{AM}_k$ of height $h(w_2)=2$.
This produces the moment 
$$
\psi(\widetilde{w}_{1})(a_2 \cdots a_{k+1})
$$ 
that can be pulled out. Here, we use $\widetilde{w}_{2}$ which means that we reduce the height of all letters in $w_2$ by one. 
This has to be done since $\psi(u)\equiv 0$ for all words $u\in \mathpzc{AM}\setminus \mathpzc{M}$, so in order to 
reproduce the moments at sites of colors $>1$ we need to reduce the 
height of $u$ to make it a reduced Motzkin word.
Let us look at sites of color $1$ associated with $a_1$ and $a_{k+2}$. 
Between $a_1(1)$ and $a_{k+2}(1)$ there are only replicas of color $>1$, among them 
$a_2(2)$ and $a_{k+1}(2)$. 
Two cases need to be considered: $i_1=i_{k+2}$ and $i_1\neq i_{k+2}$. 
If $i_1=i_{k+2}$, then in between $a_1$ and $a_{k+2}$ at site of color $1$ 
we obtain $p_1^{\perp}$, which gives $a_1p_1^{\perp}a_{k+2}=a_1a_{k+2}-a_1p_1a_{k+2}$. That is why 
$a_1$ is either connected with $a_{k+2}$ or separated from it. 
The separation produces 
$$
\varphi_{i_1}(a_1)\psi(w_2)(a_{k+2} \cdots a_{n})
$$
with a minus sign, whereas the connection gives
$$
\psi(1w_2)(a_{1}a_{k+2}\cdots a_{n}).
$$
This completes the proof of our formula.
\hfill $\blacksquare$\\

\section{Lattices $\mathcal{M}(w)$} 

We need to introduce a class of noncrossing partitions that is suitable for
computing the moments of replicas. For this purpose, we equip their blocks with colors corresponding 
to letters of Motzkin words. A relation between the colors of blocks and their nearest outer blocks 
will be of importance, as it was in the context of matricial freeness \cite{[18], [L2014]}.
In our terminology, the indices representing different algebras or subalgebras are called {\it labels}, 
whereas {\it colors} correspond to `pieces' into which variables (mainly, free random variables) 
are decomposed. In this paper, these `pieces' correspond to letters of an alphabet. 
Thus, in Section 8, where we express moments in terms of boolean cumulants, we will equip 
partitions with both gadgets, colors and labels, but for the moment we just look at colors.

A partition of the set $[n]=\{1,2, \ldots, n\}$ is the collection 
$\pi_0=\{V_1,\ldots , V_k\}$ of non-empty subsets of $[n]$, such that $V_i\cap V_j=\emptyset$
if $i\neq j$ and $[n]=V_1\cup \ldots \cup V_k$. We will adopt the convenient convention 
that 
$$
V_1<\cdots <V_k,
$$
by which we understand that the smallest number in $V_j$ is smaller than the smallest number 
in $V_{j+1}$ for any $j=1, \ldots ,k-1$. It is well known 
that it is a lattice for any $n\in {\mathbb N}$ with the partial order given by the reversed refinement, namely
$\pi_0\leq \rho_0$ if and only if $\pi_0$ is a (not necessarily proper) refinement of $\rho_0$. 

We say that $\pi_0$ is {\it noncrossing} if there do 
not exist two blocks $V_i,V_j$, where $i\neq j$, such that $k,l\in V_i$, $p,q\in V_j$ 
and $k<p<l<q$. The set of noncrossing partitions of $[n]$ will be denoted ${\rm NC}(n)$. 
It is known that ${\rm NC}(n)$ is a finite lattice for any $n\in {\mathbb N}$, although is is not a sublattice of $P(n)$. 
The partition of $[n]$ consisting of one subset is usually denoted $\hat{1}_{n}$ (it is the {\it greatest element} of 
the lattice ${\rm NC}(n)$), whereas that consisting of $n$ one-element subsets (called singletons) is denoted by $\hat{0}_{n}$(it is the {\it least element} of ${\rm NC}(n)$).

Let us suppose that $\pi_0\in {\rm NC}(n)$. Then its block $V_i$ is an {\it inner} block of $V_j$ if $p<k<q$
for any $k\in V_i$ and some $p,q\in V_j$. In that case, $V_j$ is an {\it outer} block of $V_i$.
If $V_j$ does not have any outer blocks, it is called a {\it covering block}.
If $V_i$ is not a covering block, then $V_j$ is its {\it nearest outer block} iff
$V_j$ is an outer block of $V_i$ and there is no $V_k$ such that $V_k$ is an outer block of $V_i$
and $V_j$ is an outer block of $V_k$. In other words, the nearest outer block of $V_i$ is that outer block 
$V_j$ which is the closest block to $V_i$ from the top. We use the terminology introduced in \cite{[15]}. Using 
the tree-like language of \cite{[6]}, the nearest outer block of $V_i$ is its {\it parent}.
Since the nearest outer block, if it exists, is unique,
we can write in this case 
$$
V_j=o(V_i)
$$
and we will also say that $V_j$ {\it covers} $V_i$.
A subpartition of $\pi$ is called {\it inner} with respect to $V_j$ if all its blocks are inner
blocks of $V_j$. Next, we say that $V_j$ has {\it depth} $d\in {\mathbb N}$
iff there exists a sequence of nested blocks $(V_{k_1},\ldots , V_{k_{d-1}})$, such that 
$$
(V_{k_{1}},V_{k_{2}}, \ldots , V_{k_{d-1}})=(o(V_j), o(V_{k_1}), \ldots , o(V_{k_{d-2}}))
$$
and $V_{k_{d-1}}$ is a covering block. In that case we write $d(V_i)=d$. 
Thus, the covering block of $\pi$ has depth $1$, its closest inner blocks have depth 2, etc.  
We will say that $\pi_0$ is of depth $d$ if it has a block of depth $d$ and has no blocks of depth $d+1$, we will write $d(\pi_0)=d$.
A partition $\pi_0\in {\rm NC}(n)$ is {\it irreducible} iff it has exactly one covering block.
The set of noncrossing irreducible partitions of $[n]$ is denoted ${\rm NC}_{{\rm irr}}(n)$.

\begin{figure}
\unitlength=1mm
\special{em.linewidth 0.5pt}
\linethickness{0.5pt}
\begin{picture}(140.00,30.00)(8.00,105.00)

\put(26.00,125.00){$\pi$}
\put(6.00,120.00){\line(1,0){40.00}}
\put(10.00,117.00){\line(1,0){16.00}}
\put(14.00,114.00){\line(1,0){8.00}}
\put(30.00,117.00){\line(1,0){12.00}}
\put(6.00,110.00){\line(0,1){10.00}}
\put(10.00,110.00){\line(0,1){7.00}}
\put(14.00,110.00){\line(0,1){4.00}}
\put(18.00,110.00){\line(0,1){4.00}}
\put(22.00,110.00){\line(0,1){4.00}}
\put(26.00,110.00){\line(0,1){7.00}}
\put(30.00,110.00){\line(0,1){7.00}}
\put(42.00,110.00){\line(0,1){7.00}}
\put(46.00,110.00){\line(0,1){10.00}}

\put(6.00,110.00){\circle*{1.00}}
\put(10.00,110.00){\circle*{1.30}}
\put(14.00,110.00){\circle*{1.60}}
\put(18.00,110.00){\circle*{1.60}}
\put(22.00,110.00){\circle*{1.60}}
\put(26.00,110.00){\circle*{1.30}}
\put(30.00,110.00){\circle*{1.30}}
\put(34.00,110.00){\circle*{1.60}}
\put(38.00,110.00){\circle*{1.60}}
\put(42.00,110.00){\circle*{1.30}}
\put(46.00,110.00){\circle*{1.00}}

\put(5.00,106.00){$\scriptstyle{1}$}
\put(9.00,106.00){$\scriptstyle{2}$}
\put(13.00,106.00){$\scriptstyle{3}$}
\put(17.00,106.00){$\scriptstyle{3}$}
\put(21.00,106.00){$\scriptstyle{3}$}
\put(25.00,106.00){$\scriptstyle{2}$}
\put(29.00,106.00){$\scriptstyle{2}$}
\put(33.00,106.00){$\scriptstyle{3}$}
\put(37.00,106.00){$\scriptstyle{3}$}
\put(41.00,106.00){$\scriptstyle{2}$}
\put(45.00,106.00){$\scriptstyle{1}$}

\put(76.00,125.00){$\rho$}
\put(56.00,120.00){\line(1,0){40.00}}
\put(60.00,117.00){\line(1,0){32.00}}
\put(68.00,114.00){\line(1,0){4.00}}
\put(84.00,114.00){\line(1,0){4.00}}
\put(56.00,110.00){\line(0,1){10.00}}
\put(60.00,110.00){\line(0,1){7.00}}
\put(68.00,110.00){\line(0,1){4.00}}
\put(72.00,110.00){\line(0,1){4.00}}
\put(76.00,110.00){\line(0,1){7.00}}
\put(80.00,110.00){\line(0,1){7.00}}
\put(84.00,110.00){\line(0,1){4.00}}
\put(88.00,110.00){\line(0,1){4.00}}
\put(92.00,110.00){\line(0,1){7.00}}
\put(96.00,110.00){\line(0,1){10.00}}

\put(56.00,110.00){\circle*{1.00}}
\put(60.00,110.00){\circle*{1.30}}
\put(64.00,110.00){\circle*{1.60}}
\put(68.00,110.00){\circle*{1.60}}
\put(72.00,110.00){\circle*{1.60}}
\put(76.00,110.00){\circle*{1.30}}
\put(80.00,110.00){\circle*{1.30}}
\put(84.00,110.00){\circle*{1.60}}
\put(88.00,110.00){\circle*{1.60}}
\put(92.00,110.00){\circle*{1.30}}
\put(96.00,110.00){\circle*{1.00}}

\put(55.00,106.00){$\scriptstyle{1}$}
\put(59.00,106.00){$\scriptstyle{2}$}
\put(63.00,106.00){$\scriptstyle{3}$}
\put(67.00,106.00){$\scriptstyle{3}$}
\put(71.00,106.00){$\scriptstyle{3}$}
\put(75.00,106.00){$\scriptstyle{2}$}
\put(79.00,106.00){$\scriptstyle{2}$}
\put(83.00,106.00){$\scriptstyle{3}$}
\put(87.00,106.00){$\scriptstyle{3}$}
\put(91.00,106.00){$\scriptstyle{2}$}
\put(95.00,106.00){$\scriptstyle{1}$}

\put(120.00,125.00){$\pi\vee \rho$}
\put(106.00,120.00){\line(1,0){40.00}}
\put(110.00,117.00){\line(1,0){32.00}}
\put(114.00,114.00){\line(1,0){8.00}}
\put(134.00,114.00){\line(1,0){4.00}}
\put(106.00,110.00){\line(0,1){10.00}}
\put(110.00,110.00){\line(0,1){7.00}}
\put(114.00,110.00){\line(0,1){4.00}}
\put(118.00,110.00){\line(0,1){4.00}}
\put(122.00,110.00){\line(0,1){4.00}}
\put(126.00,110.00){\line(0,1){7.00}}
\put(130.00,110.00){\line(0,1){7.00}}
\put(134.00,110.00){\line(0,1){4.00}}
\put(138.00,110.00){\line(0,1){4.00}}
\put(142.00,110.00){\line(0,1){7.00}}
\put(146.00,110.00){\line(0,1){10.00}}

\put(106.00,110.00){\circle*{1.00}}
\put(110.00,110.00){\circle*{1.30}}
\put(114.00,110.00){\circle*{1.60}}
\put(118.00,110.00){\circle*{1.60}}
\put(122.00,110.00){\circle*{1.60}}
\put(126.00,110.00){\circle*{1.30}}
\put(130.00,110.00){\circle*{1.30}}
\put(134.00,110.00){\circle*{1.60}}
\put(138.00,110.00){\circle*{1.60}}
\put(142.00,110.00){\circle*{1.30}}
\put(146.00,110.00){\circle*{1.00}}

\put(105.00,106.00){$\scriptstyle{1}$}
\put(109.00,106.00){$\scriptstyle{2}$}
\put(113.00,106.00){$\scriptstyle{3}$}
\put(117.00,106.00){$\scriptstyle{3}$}
\put(121.00,106.00){$\scriptstyle{3}$}
\put(125.00,106.00){$\scriptstyle{2}$}
\put(129.00,106.00){$\scriptstyle{2}$}
\put(133.00,106.00){$\scriptstyle{3}$}
\put(137.00,106.00){$\scriptstyle{3}$}
\put(141.00,106.00){$\scriptstyle{2}$}
\put(145.00,106.00){$\scriptstyle{1}$}

\end{picture}
\caption{Example of join in $\mathcal{M}(w)$.}
\end{figure}
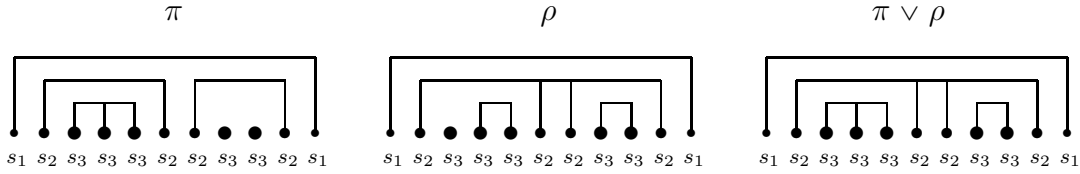

\begin{Example}
{\rm Consider the noncrossing partition $\pi_0$ in Fig.2. For the moment, let us ignore
the letters assigned to the legs of $\pi_0$ and let us treat this partition 
as an element of ${\rm NC}(11)$. We have
$V_{1}=\{1,11\}$, $V_{2}=\{2,6\}$, $V_{3}=\{3,4,5\}$,
$V_{4}=\{7,10\}$, $V_{5}=\{8\}$, $V_{6}=\{9\}$. 
The block $V_1$ is the only covering block and
$$
V_1=o(V_2)=o(V_4), \;V_2=o(V_3), \,V_4=o(V_5)=o(V_6), 
$$
with
$$
d(V_1)=1,\; d(V_2)=d(V_4)=2, \;d(V_3)=d(V_5)=d(V_6)=3.
$$
Therefore, $d(\pi_0)=3$. Note also that $\pi_0$ is irreducible. 
}
\end{Example}

We are ready to implement the structure of Motzkin paths (or, Motzkin words) into the framework 
of noncrossing partitions. The main idea is to define a subset of ${\rm NC}(n)$ for each Motzkin word $w$ of lenght $n$ that 
would consist of partitions which are adapted to $w$ in a suitable sense.

\begin{Definition}
{\rm Let $w\in \mathpzc{M}$ be of the form $w=j_1\cdots j_n$, 
where $j_1, \ldots , j_n\in {\mathbb N}$. By a {\it noncrossing partition monotonically adapted to} $w$ 
we shall understand the set of pairs 
$$
\pi=\{(V_1,v_1), \ldots, (V_k,v_k)\},
$$
where $\pi_{0}:=\{V_{1},\ldots , V_{k}\}\in {\rm NC}(n)$ and $v_j=d^{|V_j|}$ for any $j$, with $d=d(V_j)$.
The set of noncrossing partitions which are 
monotonically adapted to $w$ will be denoted by $\mathcal{M}(w)$. 
By $\mathcal{M}_{{\rm irr}}(w)$ we denote its subset consisting of irreducible partitions.
}
\end{Definition}

\begin{Definition}
{\rm Let $\pi\in \mathcal{M}(w)$ as above, which we identify with 
$$
\pi=\{v_1,\ldots , v_k\},
$$
where we assume that the corresponding blocks $V_j$ satisfy $V_1<\cdots <V_k$. Then
\begin{enumerate}[(a)]
\item
each word $v_j$ will be called a {\it block} of $\pi$,
\item
the number $d(v_j):=d(V_j)$ will be called the {\it depth} of $v_j$,
\item
if $V_j=o(V_i)$, we will write $v_j=o(v_i)$ and $(V_j,v_j)=o(V_i,v_i)$ and $v_j=o(v_i)$ and we will
say that $v_j$ is the {\it nearest outer block} of $v_i$, or that $v_j$ {\it covers} $v_i$.
\end{enumerate}
}
\end{Definition}

\begin{Example}
{\rm 
In Fig.~3, we assigned letters from the alphabet ${\mathbb N}$ to all legs of the considered partitions.
For instance, the words corresponding to the blocks of $\pi$ in Fig.~3 are: $v_1=1^2$, $v_2=2^2$, $v_3=3^3$, 
$v_4=2^2$, $v_5=3$, $v_6=3$. In fact, all partitions in Fig.~3 are monotonically adapted to $w$ since 
all words corresponding to blocks are constant. 
However, we could assign different words to the same diagrams and then 
the corresponding partitions would be adapted, but not monotonically adapted to $w$.
We choose the convention that the size of each circle is `proportional' to the height 
of the corresponding letter. This will be helpful later, when 
we assign both labels and colors to legs.
}
\end{Example}

\begin{figure}
\unitlength=1mm
\special{em:linewidth 1pt}
\linethickness{0.5pt}
\begin{picture}(180.00,110.00)(-30.00,-15.00)

\put(-10.00,85.00){$w$}
\put(26.00,85.00){$\mathcal{M}_{{\rm irr}}(w)$}
\put(100.00,85.00){$c(w)$}

\put(-10.00,70.50){$1^5$}
\put(26.00,70.00){\circle*{1.00}}
\put(29.00,70.00){\circle*{1.00}}
\put(32.00,70.00){\circle*{1.00}}
\put(35.00,70.00){\circle*{1.00}}
\put(38.00,70.00){\circle*{1.00}}
\put(26.00,70.00){\line(0,1){5}}
\put(29.00,70.00){\line(0,1){5}}
\put(32.00,70.00){\line(0,1){5}}
\put(35.00,70.00){\line(0,1){5}}
\put(38.00,70.00){\line(0,1){5}}
\put(26.00,75.00){\line(1,0){12}}
\put(100.00,71.00){$1$}

\put(-10.00,60.50){$1^321$}
\put(26.00,60.00){\circle*{1.00}}
\put(29.00,60.00){\circle*{1.00}}
\put(32.00,60.00){\circle*{1.00}}
\put(35.00,60.00){\circle*{1.30}}
\put(38.00,60.00){\circle*{1.00}}
\put(26.00,60.00){\line(0,1){5}}
\put(29.00,60.00){\line(0,1){5}}
\put(32.00,60.00){\line(0,1){5}}
\put(38.00,60.00){\line(0,1){5}}
\put(26.00,65.00){\line(1,0){12}}
\put(100.00,61.00){$1$}

\put(-10.00,50.50){$1^221^2$}
\put(26.00,50.00){\circle*{1.00}}
\put(29.00,50.00){\circle*{1.00}}
\put(32.00,50.00){\circle*{1.30}}
\put(35.00,50.00){\circle*{1.00}}
\put(38.00,50.00){\circle*{1.00}}
\put(26.00,50.00){\line(0,1){5}}
\put(29.00,50.00){\line(0,1){5}}
\put(35.00,50.00){\line(0,1){5}}
\put(38.00,50.00){\line(0,1){5}}
\put(26.00,55.00){\line(1,0){12}}
\put(100.00,51.00){$1$}

\put(-10.00,40.50){$121^3$}
\put(26.00,40.00){\circle*{1.00}}
\put(29.00,40.00){\circle*{1.30}}
\put(32.00,40.00){\circle*{1.00}}
\put(35.00,40.00){\circle*{1.00}}
\put(38.00,40.00){\circle*{1.00}}
\put(26.00,40.00){\line(0,1){5}}
\put(32.00,40.00){\line(0,1){5}}
\put(35.00,40.00){\line(0,1){5}}
\put(38.00,40.00){\line(0,1){5}}
\put(26.00,45.00){\line(1,0){12}}
\put(100.00,41.00){$1$}

\put(-10.00,30.50){$12121$}
\put(26.00,30.00){\circle*{1.00}}
\put(29.00,30.00){\circle*{1.30}}
\put(32.00,30.00){\circle*{1.00}}
\put(35.00,30.00){\circle*{1.30}}
\put(38.00,30.00){\circle*{1.00}}
\put(26.00,30.00){\line(0,1){5}}
\put(32.00,30.00){\line(0,1){5}}
\put(38.00,30.00){\line(0,1){5}}
\put(26.00,35.00){\line(1,0){12}}
\put(100.00,31.00){$1$}

\put(-10.00,20.50){$1^22^21$}
\put(26.00,20.00){\circle*{1.00}}
\put(29.00,20.00){\circle*{1.00}}
\put(32.00,20.00){\circle*{1.30}}
\put(35.00,20.00){\circle*{1.30}}
\put(38.00,20.00){\circle*{1.00}}
\put(26.00,20.00){\line(0,1){5}}
\put(29.00,20.00){\line(0,1){5}}
\put(32.00,20.00){\line(0,1){3}}
\put(35.00,20.00){\line(0,1){3}}
\put(38.00,20.00){\line(0,1){5}}
\put(26.00,25.00){\line(1,0){12}}
\put(32.00,23.00){\line(1,0){3}}

\put(43.00,20.00){\circle*{1.00}}
\put(46.00,20.00){\circle*{1.00}}
\put(49.00,20.00){\circle*{1.30}}
\put(52.00,20.00){\circle*{1.30}}
\put(55.00,20.00){\circle*{1.00}}
\put(43.00,20.00){\line(0,1){5}}
\put(46.00,20.00){\line(0,1){5}}
\put(55.00,20.00){\line(0,1){5}}
\put(43.00,25.00){\line(1,0){12}}
\put(100.00,21.00){$2$}

\put(-10.00,10.50){$12^21^2$}
\put(26.00,10.00){\circle*{1.00}}
\put(29.00,10.00){\circle*{1.30}}
\put(32.00,10.00){\circle*{1.30}}
\put(35.00,10.00){\circle*{1.00}}
\put(38.00,10.00){\circle*{1.00}}
\put(26.00,10.00){\line(0,1){5}}
\put(29.00,10.00){\line(0,1){3}}
\put(32.00,10.00){\line(0,1){3}}
\put(35.00,10.00){\line(0,1){5}}
\put(38.00,10.00){\line(0,1){5}}
\put(26.00,15.00){\line(1,0){12}}
\put(29.00,13.00){\line(1,0){3}}

\put(43.00,10.00){\circle*{1.00}}
\put(46.00,10.00){\circle*{1.30}}
\put(49.00,10.00){\circle*{1.30}}
\put(52.00,10.00){\circle*{1.00}}
\put(55.00,10.00){\circle*{1.00}}
\put(43.00,10.00){\line(0,1){5}}
\put(52.00,10.00){\line(0,1){5}}
\put(55.00,10.00){\line(0,1){5}}
\put(43.00,15.00){\line(1,0){12}}
\put(100.00,11.00){$2$}

\put(-10.00,0.50){$12^31$}
\put(26.00,0.00){\circle*{1.00}}
\put(29.00,0.00){\circle*{1.30}}
\put(32.00,0.00){\circle*{1.30}}
\put(35.00,0.00){\circle*{1.30}}
\put(38.00,0.00){\circle*{1.00}}
\put(26.00,0.00){\line(0,1){5}}
\put(29.00,0.00){\line(0,1){3}}
\put(32.00,0.00){\line(0,1){3}}
\put(35.00,0.00){\line(0,1){3}}
\put(38.00,0.00){\line(0,1){5}}
\put(26.00,5.00){\line(1,0){12}}
\put(29.00,3.00){\line(1,0){6}}

\put(43.00,0.00){\circle*{1.00}}
\put(46.00,0.00){\circle*{1.30}}
\put(49.00,0.00){\circle*{1.30}}
\put(52.00,0.00){\circle*{1.30}}
\put(55.00,0.00){\circle*{1.00}}
\put(43.00,0.00){\line(0,1){5}}
\put(46.00,0.00){\line(0,1){3}}
\put(49.00,0.00){\line(0,1){3}}
\put(55.00,0.00){\line(0,1){5}}
\put(43.00,5.00){\line(1,0){12}}
\put(46.00,3.00){\line(1,0){3}}

\put(60.00,0.00){\circle*{1.00}}
\put(63.00,0.00){\circle*{1.30}}
\put(66.00,0.00){\circle*{1.30}}
\put(69.00,0.00){\circle*{1.30}}
\put(72.00,0.00){\circle*{1.00}}
\put(60.00,0.00){\line(0,1){5}}
\put(66.00,0.00){\line(0,1){3}}
\put(69.00,0.00){\line(0,1){3}}
\put(72.00,0.00){\line(0,1){5}}
\put(60.00,5.00){\line(1,0){12}}
\put(66.00,3.00){\line(1,0){3}}

\put(77.00,0.00){\circle*{1.00}}
\put(80.00,0.00){\circle*{1.30}}
\put(83.00,0.00){\circle*{1.30}}
\put(86.00,0.00){\circle*{1.30}}
\put(89.00,0.00){\circle*{1.00}}
\put(77.00,0.00){\line(0,1){5}}
\put(89.00,0.00){\line(0,1){5}}
\put(77.00,5.00){\line(1,0){12}}
\put(100.00,01.00){$4$}

\put(-10.00,-09.50){$12321$}
\put(26.00,-10.00){\circle*{1.00}}
\put(29.00,-10.00){\circle*{1.30}}
\put(32.00,-10.00){\circle*{1.60}}
\put(35.00,-10.00){\circle*{1.30}}
\put(38.00,-10.00){\circle*{1.00}}
\put(26.00,-10.00){\line(0,1){5}}
\put(29.00,-10.00){\line(0,1){3}}
\put(35.00,-10.00){\line(0,1){3}}
\put(38.00,-10.00){\line(0,1){5}}
\put(26.00,-5.00){\line(1,0){12}}
\put(29.00,-7.00){\line(1,0){6}}
\put(100.00,-9.00){$1$}
\end{picture}
\caption{Lattices $\mathcal{M}_{{\rm irr}}(w)$ for $w\in \mathpzc{M}_{5}$. 
The corresponding cardinalities give the decomposition of the Catalan number $C_4=14=\sum_{w\in \mathpzc{M}_{5}}c(w)$.}
\end{figure}

\begin{Proposition}
The sets $\mathcal{M}(w)$ are disjoint and there is a natural bijection 
$$
{\rm NC}(n)\cong \bigcup_{w\in \mathpzc{M}_{n}}\mathcal{M}(w)
$$
for any natural $n$. The same holds true for the sets
${\rm NC}_{{\rm irr}}(n)$ and $\mathcal{M}_{{\rm irr}}(w)$.
Moreover, we have the decomposition 
$$
C_{n}=\sum_{w\in \mathpzc{M}_{n+1}}c(w),
$$
where $C_{n}=\frac{1}{n+1}{2n \choose n}$ is the nth Catalan number and 
$c(w)=|\mathcal{M}_{{\rm irr}}(w)|$ for $w\in  \mathpzc{M}_{n+1}$.
\end{Proposition}
{\it Proof.}
If $w=j_1\cdots j_n\neq j_1'\cdots j_n'=w'$, then there exists at least one 
$k$ such that $j_k\neq j_k'$. If $\pi\in \mathcal{M}(w)$ and $\pi'\in \mathcal{M}(w')$,
then the block which contains $j_k$ must be of different depth than the block that contains $j_k'$.
Hence, $\mathcal{M}(w)\cap \mathcal{M}(w')=\emptyset$.
The bijection $\tau:{\rm NC}(n)\rightarrow \bigcup_{w\in \mathpzc{M}_{n}}\mathcal{M}(w)$ is also
easy to construct. In fact, it suffices to consider irreducible partitions since each noncrossing partition 
can be decomposed into irreducible components and the depths of blocks refer to
these components.
If $\pi_{0}=\{V_1, \ldots , V_p\} \in {\rm NC}_{{\rm irr}}(n)$, there is a unique 
reduced Motzkin word $w\in \mathpzc{M}_{n}$, given by $w=j_1\ldots j_n$, where
$j_i$ is the depth of the block containing $i$, and then the blocks of $\pi\in \mathcal{M}_{{\rm irr}}(w)$ 
are of the form $(V_j,v_j)$, where $v_j$ is a constant word with color equal to the depth of $V_j$. It is clear that the 
mapping $\tau$ given by $\tau(\pi_{0})=\pi$ is a bijection. 
Finally, to show the decomposition of Catalan numbers, we first use the fact that $C_n=|{\rm NC}_{{\rm irr}}(n+1)|$, which 
follows from the formula for free cumulants of a single variable 
$$
r_n=\sum_{\pi_0\in {\rm NC}{{\rm irr}}}(-1)^{|\pi_0|-1}\beta_{\pi_0}
$$
since each boolean cumulant $\beta_{\pi_0}$ produces the coefficient $\prod_{V\in \pi_0}(-1)^{|V|-1}$ 
standing by $m_1^{n}$ and thus the contribution from each summand in the above sum to the value  
$\mu(\hat{0}_n,\hat{1}_n)$ of the M\"{o}bius function for lattice ${\rm NC}(n)$ is 
$$
(-1)^{|\pi_0|-1}(-1)^{|V_1|+\cdots +|V_p|-p}=(-1)^{n-1},
$$
where $\pi_0=\{V_1, \ldots ,V_n\}$. Since we know that $\mu(0_n,1_n)=(-1)^{n-1}C_{n-1}$, 
we conclude from the bijection discussed above that the numbers $|c(w)|$ add up to $C_{n-1}$. This completes the proof.
\hfill $\blacksquare$\\

Therefore, we have a canonical decomposition of the lattice ${\rm NC}(n)$ in terms of lattices 
$\mathcal{M}(w)$. The same holds true for the corresponding subsets of irreducible partitions. 
This fact is useful in the computations of moments and boolean cumulants of free random variables. 
Let us also observe that the cardinalities of $\mathcal{M}_{{\rm irr}}(w)$ give decompositions of Catalan numbers.
Note also that $\mathcal{M}(w)$ becomes a poset with partial order induced from ${\rm NC}(|w|)$. It turns out that 
it is a lattice.

\begin{Lemma}
For any nonempty Motzkin word $w$, the poset $\mathcal{M}(w)$ 
has a least element and a greatest element.
\end{Lemma}
{\it Proof.}
Let $w$ be a nonempty reduced Motzkin word with highest letter $d$, where $d\in \mathbb{N}$.
Let us construct a partition $\hat{0}_{w}$ that will be the least element of $\mathcal{M}(w)$. 
For that purpose, we first require each letter $d$ to be a singleton in $\pi$ 
(note that any connections between these singletons would produce a coarser partition).
Now, we cover all these letters with a minimal number of 
blocks of type $(d-1)^{k}$, where $k\geq 2$. Namely, if $w$ contains a subword of the form
$$
u=(d-1)d^{j_1}(d-1)\ldots d^{j_{k}}(d-1),
$$
where $j_1, \ldots , j_{k}, k\in \mathbb{N}$, and such $k$ is maximal, then we connect all letters 
$d-1$ in one block $(d-1)^{k+1}$. We do this for each subword of the above form 
and the remaining letters $d-1$ (which do not surround powers of $d$) are made singletons. 
Thus, the depth of each singleton $d$ is bigger by one than the depth of the constructed
block $(d-1)^{k+1}$. We continue this procedure by a reverse induction with respect to
the decreasing height of letters. Thus, if $w$ contains a subword of the form 
$$
u=jw_{1}j\ldots w_{k}j,
$$ 
where $w_1, \ldots , w_{k}$ are nonempty Motzkin words of height $j+1$ 
(in which blocks of type $i^m$, where $i>j$, have already been constructed),
we connect all letters $j$ in one block $j^{k+1}$, provided 
such $k$ is maximal. Thus, the depth of each previously constructed block 
is bigger by one than the depth of this $j^{k}$. The remaining letters $j$ are 
made singletons. The procedure ends when we reach $j=1$. By construction, 
all blocks of the constructed partition are constant and the depth of any block 
associated with a subword of height $j$ is equal to $j$. Therefore, $\pi\in \mathcal{M}(w)$. Moreover, 
removing any connection between letters of $w$ would reduce the number of outer blocks 
of at least one block and thus would produce a partition which is not in $\mathcal{M}(w)$.
Therefore, our partition is a minimal element.
It remains to justify that it is smaller than any partition in $\mathcal{M}(w)$.
If $w$ is constant, the assertion is obvious. If it is not constant, then 
suppose that there is some $\pi\in \mathcal{M}(w)$, for which it is not true that $\hat{0}_{w}\leq \pi$.
If there was a subword in $\pi$ of the form
$$
u=1w_{1}1\ldots w_{k}1,
$$ 
where $w_1, \ldots, w_k$ are nonempty Motzkin words of height $2$ and such $k$ is maximal,
and the letters $1$ in $u$ are not in one block, then one of the subwords $w_j$ or a subword of $w_j$ 
would be a block of depth $1$, which is not possible since the height of each of these subwords would be 
greater than 1. Therefore, the letters $1$ in such a word must be in the same block (of course, this block can be bigger, or longer, 
than $1^{k+1}$). The same argument applies to words $w_1, \ldots , w_k$ and its subwords of the form
$$
u=2w_{1}2\ldots w_{k}2,
$$
where $w_1, \ldots , w_k$ are nonempty Motzkin words of height 3, and so on. Therefore, any $\pi\in \mathcal{M}(w)$ must 
be bigger than $\hat{0}_{w}$. The construction of the greatest element $\hat{1}_{w}$ is similar. 
Going through similar steps as above, we build one block from all letters $1$. 
Then we build one block from all letters $2$ lying between two consecutive letters $1$. We continue in this fashion
until letters of all colors are connected in blocks. The details are left to the reader.
This completes the proof.
\hfill $\blacksquare$

\begin{Theorem}
For any nonempty Motzkin word $w$, the poset $\mathcal{M}(w)$ is a lattice.
\end{Theorem}
{\it Proof.}
Since we already know that $\mathcal{M}(w)$ has a least element, it suffices to show that for any $\pi, \rho\in  \mathcal{M}(w)$
there is a least upper bound $\pi\vee \rho$. The construction of $\pi\vee \rho$ is based on taking 
the boolean join of certain one-color subpartitions of $\pi$ and $\rho$. Namely, let
$$
\pi_j=\pi|_{w_j}\;\;\;{\rm and}\;\;\; \rho_j=\rho|_{w_j}
$$
be the subpartitions of $\pi$ and $\rho$ consisting of constant words built from $j$, which can be viewed as 
restrictions of $\pi$ and $\rho$, respectively, to the subwords $w_j$ of $w$ containing all letters $j$, where $1\leq j \leq d$ 
and $d=d(\pi)=d(\rho)$ (by definition, any partition from $\mathcal{M}(w)$ must have the same depth).
Since constant blocks built from letter $j$ must have depth $j$ by the definition of $\mathcal{M}(w)$, 
$\pi_j$ and $\rho_j$ are interval partitions of $w_j$. Moreover, blocks of both $\pi_j$ and $\rho_j$ cannot 
connect letters $j$ which are separated by letter $j-1$. Therefore, the least upper bound of $\pi_j$ 
and $\rho_j$, treated as elements of ${\mathcal Int}(w_j)$, can be identfied with the boolean join 
$
\pi_j\curlyvee \rho_j
$ 
in $I(|w_j|)$.
If we treat this operation as an operation on $\pi$ and $\rho$ involving only their blocks of color $j$, 
we can see that it does not change the depths of any of their blocks and thus the partitions obtained from 
$\pi$ and $\rho$ after each such operation are also in $\mathcal{M}(w)$. It is convenient to start 
this procedure from $w_1$ and end with $w_d$. 
If we carry out this process for all colors $j$, then we set
$$
\pi\vee \rho=\{\pi_j\curlyvee\rho_j: \; 1\leq j \leq d\},
$$
where $d=d(\pi)=d(\rho)$. By construction, $\pi\vee \rho$ is noncrossing, 
the depths of all letters do not change when we take the boolean join of $\pi_j$ and $\rho_j$ 
for all $j$ and that is why $\pi\vee \rho\in \mathcal{M}(w)$. 
Moreover, it is the least upper bound of $\pi$ and $\rho$ since otherwise there would exist $j$ such that 
$\pi_{j}\curlyvee\rho_{j}$ would not be the least upper bound of $\pi_{j}$ and $\rho_j$. This would
reduce to a situation in which the restriction of $\pi_j\curlyvee \rho_j$ to each union of 
blocks of color $j$ which have the same nearest outer block, say some $u_j\subset w_j$, 
would not be the least upper bound of the restrictions of $\pi_j$ and $\rho_j$ to $u_j$. But these blocks
have the same depth and the only way to join them is to take the union of neighboring 
blocks since they have to retain the same depth. Therefore, the least upper bound of each such restriction
must by their boolean join. This completes the proof.
\hfill $\blacksquare$\\

\section{Boolean cumulant formula}

It can be seen from the considerations of Section 6 that 
one can compute effectively the moments of orthogonal replicas, using some recursions.
In particular, Lemma 6.2 indicates that the decomposition of 
Motzkin paths related to the recurrence formula for Motzkin numbers 
might lead to a combinatorial formula in which partitions from the 
lattices $\mathcal{M}(w)$ play the main role.

In the partition context, it is more convenient to use multilinear 
Motzkin functionals $\varphi(w)$ rather than the linear ones. This will
also be the case in our study of Motzkin cumulants \cite{[L2022]}.
Thus, let 
$$
{\mathcal A}:=*_{i\in I}\mathcal{A}_{i}
$$
and let $a_{k}\in \mathcal{A}_{i_k}$, where $i_1, \ldots, i_n\in I$ with the associated replicas 
$a_{k}(j)$, where $j\in \mathbb{N}$. 

In order to compute moments under $\varphi(w)$, we shall need the 
partitions which are not only monotonically adapted to $w$, but also to the label
$\ell=(i_1,  \ldots i_n)$ of the sequence of variables $(a_1, \ldots, a_n)$, where
$a_k\in {\mathcal A}_{i_k}$ for $k=1, \ldots, n$.
In order to introduce these subsets, we shall use the concept of
a {\it monotone chain of blocks} in $\pi\in \mathcal{M}(w)$, where $w\in \mathpzc{M}$, by which we understand
a sequence of blocks $(V_{j_1}, \ldots , V_{j_p})$ such that 
$$
V_{j_1}=o(V_{j_2}), \ldots V_{j_{p-1}}=o(V_{j_p})
$$
and for which the sequence of colors of the associated words $(v_{j_1}, \ldots , v_{j_p})$ is monotonically 
increasing, namely $h(v_{j_1})< \ldots <h(v_{j_p})$. We are ready to define families of partitions which 
are monotonically adapted to the label $\ell$ in a suitable sense (apart from being monotonically adapted to the word $w$).  

\begin{Definition}
{\rm 
We will say that $\pi\in \mathcal{M}(w)$  
is {\it monotonically adapted to} $\ell=(i_1, \ldots, i_n)$, where $n=|w|$, if
\begin {enumerate}[i)]
\item
$\pi$ is adapted to $\ell$, namely $i_{k}=i_l$ whenever $k,l\in V$ for all $V\in \pi$,
\item
the labels assigned to blocks alternate in monotone chains of blocks in $\pi$. 
\end{enumerate}
The set of such partitions will be denoted by $\mathcal{M}(w, \ell)$ and 
$\mathcal{M}_{{\rm irr}}(w,\ell)$ will be its subset consisting of irreducible partitions.
}
\end{Definition}

\begin{Definition}
{\rm 
For any $w\in \mathpzc{M}$, define a family
$\{\beta_{\pi}:\pi \in \mathcal{M}(w,\ell)\}$ of multilinear functionals by the multiplicative formula
$$
\beta_{\pi}[a_1, \ldots, a_n]:=\prod_{V\in \pi}\beta(V)[a_1, \ldots, a_n],
$$
where $a_k\in {\mathcal A}_{i_k}$ and $k=1, \ldots , n$, where $\beta(V)$ is given by Definition 2.1.
The functionals $\beta_{\pi}$ remind partitioned boolean functionals associated 
with the boolean product of states except that $\pi$ does not have to be an interval partition. 
In the formula for $\beta(V)$ the same symbol $\beta_k$ is used for the mixed boolean cumulants of 
order $k$ of variables from $({\mathcal A}_1,\varphi_1)$ or from $({\mathcal A}_{2},\varphi_2)$. Also, we shall
use the notation $\beta_k$ for mixed boolean cumulants corresponding to $\varphi_1\star \varphi_2$.
}
\end{Definition}

\begin{Example}
{\rm 
Before we state a general result, let us give some direct computations of mixed moments of orthogonal replicas
under the state $\Phi$ of Definition 3.1. 
On this occasion, a few examples of $\beta_{\pi}$ will also be given.
Let $w=121^22^21$, take the labeling $\ell=(1,2,1,1,2,2,1)$ and let 
$$
z(w)=a_1(1)b_1(2)a_3(1)a_4(1)b_5(2)b_6(2)a_7(1).
$$
We obtain
\begin{eqnarray*}
\Phi(z(w))&=&
\widetilde{\varphi}_{1}(a_1p^{\perp}_1a_3a_4p^{\perp}_1a_7)\varphi_{2}(b_2)\varphi_{2}(b_5b_6)\\
&=& \beta_3(a_1,a_3a_4,a_7)\beta_1(b_2)(\beta_2(b_5,b_6)+\beta_1(b_5)\beta_1(b_6))\\
&=&
(\beta_4(a_1,a_3,a_4,a_7)+\beta_2(a_1,a_3)\beta_2(a_4,a_7))\\
&\times& \beta_1(b_2)(\beta_2(b_5,b_6)+\beta_1(b_5)\beta_1(b_6))\\
&=&(\beta_{\pi_1}+\beta_{\pi_2}+\beta_{\pi_3}+\beta_{\pi_4})[a_1,b_2,a_3, a_4, b_5,b_6, a_7],
\end{eqnarray*}
using Proposition 3.1 and Corollary 3.1, where the partitions are of the form\\
\unitlength=1mm
\special{em:linewidth 1pt}
\linethickness{0.5pt}
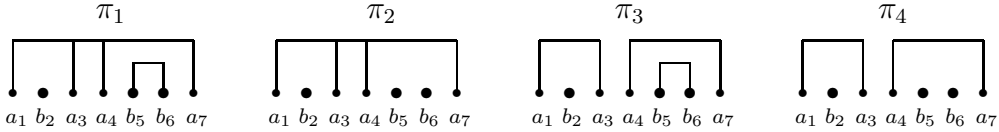
\begin{picture}(120.00,25.00)(00.00,0.00)

\put(21.00,20.00){$\pi_1$}
\put(10.00,10.00){\line(0,1){7}}
\put(18.00,10.00){\line(0,1){7}}
\put(22.00,10.00){\line(0,1){7}}
\put(26.00,10.00){\line(0,1){4}}
\put(30.00,10.00){\line(0,1){4}}
\put(34.00,10.00){\line(0,1){7}}

\put(10.00,17.00){\line(1,0){24}}
\put(26.00,14.00){\line(1,0){4}}
\put(10.00,10.00){\circle*{1.00}}
\put(14.00,10.00){\circle*{1.50}}
\put(18.00,10.00){\circle*{1.00}}
\put(22.00,10.00){\circle*{1.00}}
\put(26.00,10.00){\circle*{1.50}}
\put(30.00,10.00){\circle*{1.50}}
\put(34.00,10.00){\circle*{1.00}}

\put(9.00,06.00){$\scriptstyle{a_1}$}
\put(13.00,06.00){$\scriptstyle{b_2}$}
\put(17,06.00){$\scriptstyle{a_3}$}
\put(21.00,06.00){$\scriptstyle{a_4}$}
\put(25.00,06.00){$\scriptstyle{b_5}$}
\put(29.00,06.00){$\scriptstyle{b_6}$}
\put(33.00,06.00){$\scriptstyle{a_7}$}


\put(57.00,20.00){$\pi_2$}
\put(45.00,10.00){\line(0,1){7}}
\put(53.00,10.00){\line(0,1){7}}
\put(57.00,10.00){\line(0,1){7}}
\put(69.00,10.00){\line(0,1){7}}

\put(45.00,17.00){\line(1,0){24}}
\put(45.00,10.00){\circle*{1.00}}
\put(49.00,10.00){\circle*{1.50}}
\put(53.00,10.00){\circle*{1.00}}
\put(57.00,10.00){\circle*{1.00}}
\put(61.00,10.00){\circle*{1.50}}
\put(65.00,10.00){\circle*{1.50}}
\put(69.00,10.00){\circle*{1.00}}

\put(44.00,06.00){$\scriptstyle{a_1}$}
\put(48.00,06.00){$\scriptstyle{b_2}$}
\put(52.00,06.00){$\scriptstyle{a_3}$}
\put(56.00,06.00){$\scriptstyle{a_4}$}
\put(60.00,06.00){$\scriptstyle{b_5}$}
\put(64.00,06.00){$\scriptstyle{b_6}$}
\put(68.00,06.00){$\scriptstyle{a_7}$}


\put(90.00,20.00){$\pi_3$}
\put(80.00,10.00){\line(0,1){7}}
\put(88.00,10.00){\line(0,1){7}}
\put(92.00,10.00){\line(0,1){7}}
\put(96.00,10.00){\line(0,1){4}}
\put(100.00,10.00){\line(0,1){4}}
\put(104.00,10.00){\line(0,1){7}}

\put(80.00,17.00){\line(1,0){8}}
\put(92.00,17.00){\line(1,0){12}}
\put(96.00,14.00){\line(1,0){4}}
\put(80.00,10.00){\circle*{1.00}}
\put(84.00,10.00){\circle*{1.50}}
\put(88.00,10.00){\circle*{1.00}}
\put(92.00,10.00){\circle*{1.00}}
\put(96.00,10.00){\circle*{1.50}}
\put(100.00,10.00){\circle*{1.50}}
\put(104.00,10.00){\circle*{1.00}}

\put(79.00,06.00){$\scriptstyle{a_1}$}
\put(83.00,06.00){$\scriptstyle{b_2}$}
\put(87.00,06.00){$\scriptstyle{a_3}$}
\put(91.00,06.00){$\scriptstyle{a_4}$}
\put(95.00,06.00){$\scriptstyle{b_5}$}
\put(99.00,06.00){$\scriptstyle{b_6}$}
\put(103.00,06.00){$\scriptstyle{a_7}$}


\put(125.00,20.00){$\pi_4$}
\put(115.00,10.00){\line(0,1){7}}
\put(123.00,10.00){\line(0,1){7}}
\put(127.00,10.00){\line(0,1){7}}
\put(139.00,10.00){\line(0,1){7}}

\put(115.00,17.00){\line(1,0){8}}
\put(127.00,17.00){\line(1,0){12}}
\put(115.00,10.00){\circle*{1.00}}
\put(119.00,10.00){\circle*{1.50}}
\put(123.00,10.00){\circle*{1.00}}
\put(127.00,10.00){\circle*{1.00}}
\put(131.00,10.00){\circle*{1.50}}
\put(135.00,10.00){\circle*{1.50}}
\put(139.00,10.00){\circle*{1.00}}

\put(114.00,06.00){$\scriptstyle{a_1}$}
\put(118.00,06.00){$\scriptstyle{b_2}$}
\put(122.00,06.00){$\scriptstyle{a_3}$}
\put(126.00,06.00){$\scriptstyle{a_4}$}
\put(130.00,06.00){$\scriptstyle{b_5}$}
\put(134.00,06.00){$\scriptstyle{b_6}$}
\put(138.00,06.00){$\scriptstyle{a_7}$}

\end{picture}
\\
with the letters of $w$ omitted for simplicity.
It can be seen that 
$$
\mathcal{M}(w,\ell)=\{\pi_1,\pi_2, \pi_3, \pi_4\}
$$
which will agree with the formula of Theorem 8.1 given below.
Let us observe that the same diagrams can be associated with different labelings.
Thus, if we are given a noncrossing partition $\pi$ of $w$, we can look for labelings 
to which $\pi$ is adapted and which alternate in monotone chains of blocks,
denoted $\mathcal{L}(\pi)$. It can be seen from the above that we have the following 
cardinalities of the sets of possible labelings:
$|\mathcal{L}(\pi_1)|=2$, $|\mathcal{L}(\pi_2)|=2$, $|\mathcal{L}(\pi_3)|=4$, $|\mathcal{L}(\pi_4)|=4$.
}
\end{Example}

We have expressed the moment $\Phi(z(w))$ in term of boolean cumulants. This is quite natural in view of 
the role played by projections $p_j^{\perp}$ in computations involving 
these cumulants as shown in Section 3. At the same time, it can be seen that 
this `orthogonality effect' plays the main role in the recurrence of Lemma 6.2.
We would like to bridge together these two observations.

\begin{Theorem}
For any $w=j_1\ldots j_n\in \mathpzc{M}_{n}$ and $n\in {\mathbb N}$, 
let $a_k\in {\mathcal A}_{i_k}$ for any $k=1, \ldots , n$, where $i_1, \ldots, i_n\in I$.
It holds that
\begin{eqnarray*}
\varphi(w)(a_1, \ldots, a_n)
&=&
\sum_{\pi\in \mathcal{M}(w, \ell)}\beta_{\pi}[a_1, \ldots , a_n],
\end{eqnarray*}
where $\ell=\ell(a_1, \ldots , a_n)$.
\end{Theorem}
{\it Proof.}
Since $w\in \mathpzc{M}_n$, it can be written in the form
$$
w=1^{n_0}w_11^{n_1}\cdots w_{p}1^{n_p},
$$
where $n_0,n_1, \ldots , n_p\in {\mathbb N}$, with $n_0+\cdots +n_p+|w_1|+\cdots +|w_p|=n$, and 
subwords $w_1, \ldots , w_p\in \mathpzc{AM}$ are of height $2$.
We will use Proposition 3.2 and Lemma 6.2 to compute $\varphi(w)(a_1, \ldots, a_n)$. 
Let us point out that we do not assume that $i_1\neq \cdots \neq i_n$.
Denote the corresponding product of replicas by
$$
z(w)=a_{1}(j_1)\cdots a_{n}(j_n)
$$
and let $z(u)$ be the product of replicas corresponding to any subword $u$ of $w$. Thus,
$$
z(w)=z(1^{n_0})z(w_1)z(1^{n_1})\cdots z(w_p)z(1^{n_p}).
$$ 
By Proposition 3.2, 
$$
\varphi(w)(a_1, \ldots , a_n)=\Phi(z(w))=0
$$ 
unless to each subword $1w_k1$ we assign $a(1)z(w_k)a'(1)$, where $a,a'$ are from the same 
algebra and $z(w_k)$ does not have a replica $a''(2)$ for $a''$ from the same algebra.
Moreover, the neighboring replicas from the same algebra must be of the same color.
In all these cases $\mathcal{M}(w,\ell)= \emptyset$, which means that our statement holds.
If the neighboring replicas are from the same algebra and have the same color, i.e.
$j_i=j_{i+1}$ and $i_{k}=i_{k+1}$, then we can use the homomorphic property of Lemma 5.1
to compress the word $w$ to a word $w'$ of lenght $m<n$ for which $i_1\neq \cdots \neq i_m$.
Then we apply Lemma 6.2 to $w'$. Of course, these compressions take place either within 
some $1^{n_q}$ or within some $w_q$. Then we use Lemma 6.1 which allows us to factor out a 
product of moments of replicas associated with $1^{n_0}$ which have 
alternating labels (if there is one label, we just compress $1^{n_0}$ to $1$).
Supposing that no compressions are needed for $1^{n_0}$ since the neighboring replicas have different labels,
we arrive at
$$
\prod_{k=1}^{n_{0}-1}\varphi_{i_k}(a_k)\psi(1w_11^{n_1}\cdots w_p 1^{n_p})(a_{n_{0}}\cdots a_n)
$$
to which we can apply Lemma 6.2 to pull out the moment associated with $w_1'$, a compression of $w_1$.
The same is done next for $1^{n_1}, w_2, \ldots, w_p, 1^{n_p}$.  
At the end we decompress all subwords to reproduce the arguments of $\varphi(w)$. 
In effect, after some compressions, followed by a repeated application of Lemma 6.2, and 
finally some decompressions, we can pull out the product of moments of the form
$$
\prod_{k=1}^{p}\varphi(\widetilde{w}_k)(A_k),
$$
where $A_k=(a_{q}, \ldots , a_{r})$ for some $1<q<r<n$. Let us recall that $\widetilde{w}_k$ 
is the reduced word obtained from $w_k$ by mapping each letter $j$ onto $j-1$.
Let us now analyze what is happening at sites of color $1$ corresponding
to letter $1$. Note that Proposition 3.2 and Lemma 6.2 show that 
variables associated with letters $1$ which surround $w_1$ are either connected 
or separated (in that case a minus sign is produced).
A repeated application of Lemma 6.2 to all subsequent subwords $w_k$ gives the same effect:
the variables corresponding to the letters $1$ that surround $w_k$ 
experience the same effect: they are either connected or separated (each separation gives a minus sign).
For each $1\leq k \leq p$, these variables must belong to the same algebra.
If it is ${\mathcal A}_{i}$, it is convenient to use the projection 
$p_i^{\perp}=1-p_{i}$ to encode this fact. 
Therefore, we need to multiply the above product of moments by alternating products of moments of the form
$$
\widetilde{\varphi}_{i}(\cdots ap_{i}^{\perp}a'\cdots ),
$$
where $a,a'\in \mathcal{A}_{i}$ for suitable $i$. Now, the number of $p_{i}^{\perp}$-moments of the above type 
depends on how many times the labels of replicas of color $1$ change 
when going from left to right. Of course, $n_0+\cdots +n_p$ is equal to the total number of variables 
which appear in these $p_{i}^{\perp}$-moments. Using Proposition 3.1 and Corollary 3.1, we get
$$
\varphi(w)(a_1, \ldots, a_n)=
\sum_{\stackrel{\pi'\in {\mathcal Int}(w', \ell')}
{\scriptscriptstyle 1^{n_0}\sim \cdots \sim 1^{n_{p}}
}}
\beta_{\pi'}[a_1,\ldots , a_n]
\prod_{k=1}^{p}\varphi(\widetilde{w}_{k})(A_k),
$$
where ${\mathcal Int}(w', \ell')$ is the set of interval partitions of the word 
$$
w'=1^{n_0}\cdots 1^{n_{p}}
$$
to which the associated label $\ell'$ is adapted in the sense that 
the labels corresponding to the same block of $\pi$ are equal.
By Corollary 3.1, these partitions must connect subwords $1^{n_0}, \ldots , 1^{n_{p}}$, which we denote by
$1^{n_0}\sim \cdots \sim 1^{n_{p}}$ (notation is similar to that in Corollary 3.1). Moreover,
$$
\beta_{\pi'}[a_1, \ldots, a_n]=\prod_{v\in \pi'}\beta(v)[a_1, \ldots , a_n]
$$ 
for $\pi'\in {\mathcal Int}(w', \ell')$ ($\pi'$ will be a subpartition of a partition of $w$, hence the notation on the RHS 
is natural). It remains to justify that replicas of color $1$ cannot have the same labels as their 
nearest outer block. This follows from the conditional orthogonality of replicas (Proposition 3.2).
Therefore, we conclude that we have constructed a class of partitions $\pi$ of 
$w$ by constructing all interval partitions of $w'$ which connect subwords 
$1^{n_0}, \ldots , 1^{n_{p}}$ and completing them with 
their inner blocks built from partitions corresponding to all subwords $w_1, \ldots , w_p$ separately. 
The latter, by an induction argument, are obtained in the formulas
$$
\varphi(\widetilde{w}_{k})(A_k)
=\sum_{\pi_k\in \mathcal{M}(\widetilde{w}_{k},\widetilde{\ell}_{k})}\beta_{\pi_k}[a_{1}, \ldots ,a_{n}]
$$ 
for any $k=1, \ldots , p$, where $\widetilde{\ell}_{k}$ is the restriction of $\ell$ to 
the indices associated with $\widetilde{w}_{k}$.  In view of the arguments on labels presented above,
the partitions $\pi$ must be noncrossing and adapted to $\ell$. Moreover, the labels of replicas must change as 
the depth of blocks of $\pi$ increases. 
Thus, we have 
$$
\pi:=\pi'\cup \pi_1\cup \ldots \cup \pi_p\in \mathcal{M}(w,\ell).
$$
Of course, with the above notations, 
$$
\beta_{\pi'}[a_1, \ldots , a_n] \cdot \prod_{j=1}^{p}\beta_{\pi_j}[a_1, \ldots , a_n]=\beta_{\pi}[a_1,\ldots , a_n]
$$
and thus we obtain the desired form of the product of cumulants corresponding to $\pi$. 
It remains to justify that all partitions from $\mathcal{M}(w,\ell)$ are obtained. 
However, by induction, it suffices to justify that all interval partitions $\pi'$ of $w'$ 
which are adapted to $\ell'$ are obtained. This follows from Corollary 3.1 and the fact that 
replicas of type $a_i(1)$ or $b_j(1)$, respectively,
produce $p^{\perp}$ in between variables $b_j$ and $a_i$. Therefore, the proof is completed.
\hfill $\blacksquare$\\

\begin{Example}
{\rm 
Let $w=1232^212^21$ and consider two labelings, say $\ell, \ell'$, shown in that order, 
in the diagrams of slightly bigger depth than those in Example 8.1.
\\
\unitlength=1mm
\special{em:linewidth 1pt}
\linethickness{0.5pt}
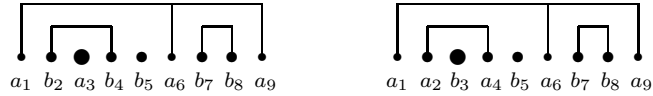
\begin{picture}(120.00,25.00)(-27.00,0.00)

\put(5.00,10.00){\line(0,1){7}}
\put(9.00,10.00){\line(0,1){4}}
\put(17.00,10.00){\line(0,1){4}}
\put(25.00,10.00){\line(0,1){7}}
\put(29.00,10.00){\line(0,1){4}}
\put(33.00,10.00){\line(0,1){4}}
\put(37.00,10.00){\line(0,1){7}}

\put(5.00,17.00){\line(1,0){32}}
\put(9.00,14.00){\line(1,0){8}}
\put(29.00,14.00){\line(1,0){4}}

\put(5.00,10.00){\circle*{1.00}}
\put(9.00,10.00){\circle*{1.50}}
\put(13.00,10.00){\circle*{2.00}}
\put(17.00,10.00){\circle*{1.50}}
\put(21.00,10.00){\circle*{1.50}}
\put(25.00,10.00){\circle*{1.00}}
\put(29.00,10.00){\circle*{1.50}}
\put(33.00,10.00){\circle*{1.50}}
\put(37.00,10.00){\circle*{1.00}}

\put(3.50,06.00){$\scriptstyle{a_1}$}
\put(8.00,06.00){$\scriptstyle{b_2}$}
\put(12.00,06.00){$\scriptstyle{a_3}$}
\put(16.06,06.00){$\scriptstyle{b_4}$}
\put(20.00,06.00){$\scriptstyle{b_5}$}
\put(24.00,06.00){$\scriptstyle{a_6}$}
\put(28,06.00){$\scriptstyle{b_7}$}
\put(32.00,06.00){$\scriptstyle{b_8}$}
\put(36.00,06.00){$\scriptstyle{a_9}$}

\put(55.00,10.00){\line(0,1){7}}
\put(59.00,10.00){\line(0,1){4}}
\put(67.00,10.00){\line(0,1){4}}
\put(75.00,10.00){\line(0,1){7}}
\put(79.00,10.00){\line(0,1){4}}
\put(83.00,10.00){\line(0,1){4}}
\put(87.00,10.00){\line(0,1){7}}

\put(55.00,17.00){\line(1,0){32}}
\put(59.00,14.00){\line(1,0){8}}
\put(79.00,14.00){\line(1,0){4}}

\put(55.00,10.00){\circle*{1.00}}
\put(59.00,10.00){\circle*{1.50}}
\put(63.00,10.00){\circle*{2.00}}
\put(67.00,10.00){\circle*{1.50}}
\put(71.00,10.00){\circle*{1.50}}
\put(75.00,10.00){\circle*{1.00}}
\put(79.00,10.00){\circle*{1.50}}
\put(83.00,10.00){\circle*{1.50}}
\put(87.00,10.00){\circle*{1.00}}

\put(53.50,06.00){$\scriptstyle{a_1}$}
\put(58.00,06.00){$\scriptstyle{a_2}$}
\put(62.00,06.00){$\scriptstyle{b_3}$}
\put(66.00,06.00){$\scriptstyle{a_4}$}
\put(70.00,06.00){$\scriptstyle{b_5}$}
\put(74.00,06.00){$\scriptstyle{a_6}$}
\put(78.00,06.00){$\scriptstyle{b_7}$}
\put(82.00,06.00){$\scriptstyle{b_8}$}
\put(86.00,06.00){$\scriptstyle{a_9}$}
\end{picture}
\\
Denote blocks as follows:
$$
V_1=\{1,6,9\}, V_2=\{2,4\}, V_3=\{3\}, V_4=\{5\}, V_5=\{7,8\}.
$$
Note that $\ell \in \mathcal{L}(\pi)$ 
since labels alternate in all monotone chains:
$(V_1,V_2,V_3)$, $(V_1,V_4)$ and $(V_1,V_5)$.
In turn, $\ell'\notin {\mathcal L}(\pi)$ since labels corresponding to blocks $V_1$ and $V_2$ in the 
monotone chain $(V_1,V_2,V_3)$ are equal and 
thus $\ell'$ is not alternating in this chain. 
If we replaced $w$ by $u=1^221212^21$ and denoted by $\rho$ 
the partition corresponding to the same $\pi_0$ (which means the same diagram), then we would 
obtain $\ell,\ell'\in \mathcal{L}(\rho)$ since in that case labels would alternate in all monotone 
chains: $(V_2, V_3)$, $(V_1,V_4)$ and $(V_1,V_5)$ (in that case, $(V_1,V_2,V_3)$ is not a monotone chain).}
\end{Example}

\begin{Definition}
{\rm Let $\ell=\ell_1\cdots \ell_n$ be a labeling, where $\ell_{i}\in I$ for any $i=1, \ldots , n$. 
We will say that $\pi_{0}=\{V_1, \ldots, V_p\}\in {\rm NC}(n)$ is {\it adapted to $\ell$} if
\begin{enumerate}[i)]
\item
labels match within blocks: $\ell_i=\ell_j$ whenever $i,j\in V_k$ for some $k$,
then we set $\ell(V_k)=\ell_{i}$ and call it the {\it label} of $V_k$,
\item
labels alternate within each chain of nested blocks: if $V_i,V_j\in \pi_0$ and $V_i=o(V_j)$, then $\ell(V_j)\neq \ell(V_i)$.
\end{enumerate}
The set of noncrossing partitions of $[n]$ adapted to $\ell$ will be denoted ${\rm NC}(n, \ell)$.
}
\end{Definition}

One of the consequences of Theorem 8.1 is a formula for moments of free random variables in terms of boolean cumulants 
similar to Theorem 1.2 derived in \cite{[6]} and Proposition 4.31 in \cite{[10]}. In other words, Theorem 8.1 can be 
viewed as a Motzkin decomposition of the formula given below.

\begin{Corollary}
With the notation of Theorem 5.1, we have the decomposition of moments of free random variables
$$
(\star_{i\in I}\varphi_{i})(x_{1} \cdots x_{n})=
\sum_{\pi_{0} \in {\rm NC}(n, \ell)}\beta_{\pi_0}[x_1,\ldots, x_n]
$$
for any $n\in {\mathbb N}$, where $x_k=\tau(a_k)$, and $a_k\in {\mathcal A}_{i_k}$ and 
$i_1, \ldots, i_n\in I$ and $k=1, \ldots, n$.
\end{Corollary}
{\it Proof.}
Recall that $\tau$ is just the unit identification mapping of Theorem 5.1. 
By Theorem 8.1, we have
\begin{eqnarray*}
(\star_{i\in I}\varphi_{i})(x_{1} \cdots  x_{n})&=&
\sum_{w\in \mathpzc{M}_{n}}
\varphi(w)(a_1, \ldots ,a_n)\\
&=&
\sum_{w\in \mathpzc{M}_{n}}
\sum_{\pi\in \mathcal{M}(w, \ell)}
\beta_{\pi}[a_1, \ldots , a_n]
\\
&=&
\sum_{\pi_{0}\in {\rm NC}(n,\ell)}
\beta_{\pi_0}[x_1,\ldots, x_n]
\end{eqnarray*}
since ${\rm NC}(n)\cong \bigcup_{w\in \mathpzc{M}_{n}}{\mathcal M}(w)$ and thus 
${\rm NC}(n,\ell)\cong \bigcup_{w\in \mathpzc{M}_{n}}{\mathcal M}(w,\ell)$, which completes the proof.
\hfill $\blacksquare$

\begin{Definition}
{\rm Let $(B_{n})$ be the sequence of boolean cumulants of orthogonal replicas
associated with the state $\Phi$. Define the family of multilinear functionals
$$
\{B(w):w\in \mathpzc{M}\},
$$
on the Cartesian powers of ${\mathcal A}$ by
$$
B(w)(a_{1}, \ldots , a_{n}):=B_n(a_1(j_1), \ldots , a_n(j_n)),
$$
where $a_k\in {\mathcal A}_{i_k}$ for $k=1, \ldots n$ and $w=j_1\cdots j_n$.}
\end{Definition}

\begin{Theorem}
The functionals $B(w)$ are given by the formula
$$
B(w)(a_1, \ldots , a_n)=\sum_{\pi\in \mathcal{M}_{{\rm irr}}(w, \ell)}\beta_{\pi}[a_1,\ldots, a_n] 
$$
where $\ell=\ell(a_1, \ldots, a_n)$. 
\end{Theorem}
{\it Proof.}
Using the moment-cumulant formula for boolean cumulants, we obtain
\begin{eqnarray*}
\varphi(w)(a_1,\ldots , a_n)&:=& \sum_{\pi_{0}\in {\rm Int}(n)}B_{\pi_{0}}[a_1(j_1), \ldots , a_n(j_n)]\\
&=& \sum_{\pi\in {\mathcal Int}(w)}B_{\pi}[a_1, \ldots , a_n]
\end{eqnarray*}
where ${\mathcal Int}(w)$ is the set of interval partitions adapted to $w$, by which we understand that 
each word $v_j\in \pi\in {\mathcal Int}(w)$ is a reduced Motzkin word. Indeed, if some $v_j$ is not a reduced 
Motzkin word, then the corresponding $B(v_j)=0$. Clearly,  $\hat{1}_{w}\in {\mathcal Int}(w)$. Thus, 
cumulants $B(w)$, which are identified with $B(\hat{1}_{w})$, are uniquely determined from the above equations
for all $w\in \mathpzc{M}$. Therefore, in view of Theorem 8.1, it suffices to show that 
$$
\varphi(w)(a_1, \ldots, a_n)=\sum_{\pi\in {\mathcal Int}(w)}\prod_{v\in \pi}\sum_{\rho_{v}\in \mathcal{M}_{{\rm irr}}(v, \ell(v))}
\beta_{\rho_v}[z_1, \ldots, z_n],
$$
where $\ell(v)$ is the label of $v$.
In this expression, we partition $w$ into consecutive reduced Motzkin subwords $v$ 
(these correspond to intervals) and then each $v$
is decomposed into irreducible partitions which are monotonically adapted to $v$ and adapted to the label $\ell(v)$, which means that 
labels of blocks alternate in monotonic chains. 
It is clear that we obtain 
a partition $\pi'\in \mathcal{M}(w, \ell)$ of the form
$$
\pi'=\bigcup_{v\in \pi}\rho_{v}.
$$
Moreover, if $\pi_1\neq \pi_2$, where $\pi_1, \pi_2\in {\mathcal Int}(w)$, 
then each $\pi_1'$ obtained in this process from $\pi_1$ is different from each $\pi_2'$ obtained from $\pi_2$ 
due to irreducibility of each $\rho_{v}\in {\mathcal M}(v, \ell(v))$  and each $\rho_{v'}\in {\mathcal M}(v', \ell(v'))$, 
for $v\in \pi_1$ and $v'\in \pi_2$.
Therefore, by decomposing different interval partitions we obtain different partitions which are 
monotonically adapted to $w$.  Now, suppose that $\pi'\in \mathcal{M}(w, \ell)$ is arbitrary. We need to show that
there is a unique $\pi\in {\mathcal Int}(n)$ and a family $\{\rho_{v}:v\in \pi\}$ which gives $\pi'$ in the decomposition process. 
If $\pi'$ is irreducible, then there is exactly one interval partition, namely $\pi=\hat{1}_{w}$ which gives
$\pi'$ in the decomposition process (there is a unique $\rho$  given by $\rho=\pi'$). 
If $\pi'$ is not irreducible, it is a union of irreducible components which correspond to intervals 
$v_1, \ldots, v_k$ such that $w=v_1\cdots v_k$ and 
$$
\pi'=\pi_1'\cup\cdots \cup \pi_k',
$$
where $\pi_j\in \mathcal{M}_{{\rm irr}}(v_j, \ell(v_j))$ for each $j$.
This defines a unique interval partition
$\pi=\{v_1, \ldots , v_k\}$ whose monotone decomposition gives $\pi'$. We treat each $\pi_j'$ in the same way as 
$\pi'$ a while ago in the case when $\pi'$ was irreducible. Next, each block corresponding to $\pi_j'$ is decomposed as a subpartition of $\pi_j'$. Call it $\rho_{v_j}$. This completes the proof. 
\hfill $\blacksquare$

\begin{Corollary}
With the above notations, we have the decomposition of boolean cumulants of free random variables
$$
\beta_n(x_{1}, \ldots , x_{n})=
\sum_{\pi_0\in {\rm NC}_{{\rm irr}}(n, \ell)}\beta_{\pi_0}[x_1,\ldots, x_n]
$$
for any $n\in {\mathbb N}$ and any $x_k\in {\mathcal A}_{i_k}$, where $k=1, \ldots, n$ and $i_1, \ldots, i_n\in I$.
\end{Corollary}
{\it Proof.}
By Theorem 5.1, multilinearity of $B$ and Theorem 8.2, and with the notation of Corollary 8.1, we have  
\begin{eqnarray*}
\beta_n(x_{1}, \ldots , x_{n})
&=&
\sum_{w\in \mathpzc{M}_{n}}
B(w)(a_1,\ldots, a_n)\\
&=&
\sum_{w\in \mathpzc{M}_{n}}
\sum_{\pi\in \mathcal{M}_{{\rm irr}}(w, \ell)}\beta_{\pi}[a_1,\ldots, a_n]\\
&=&
\sum_{\pi_{0}\in {\rm NC}_{{\rm irr}}(n, \ell)}
\beta_{\pi_{0}}[x_1,\ldots, x_n]
\end{eqnarray*}
since ${\rm NC}_{{\rm irr}}(n,\ell)\cong \bigcup_{w\in \mathpzc{M}_{n}}{\mathcal M}_{{\rm irr}}(w,\ell)$, 
which completes the proof.
\hfill $\blacksquare$

\noindent\\[10pt]
{\bf Acknowledgement}\\
I would like to express my thanks to the Reviewer for valuable comments, remarks and suggestions.

\end{document}